\documentclass[smallextended]{svjour3}

\usepackage[centertags]{amsmath}
\usepackage{amstext,amssymb,amsopn}
\usepackage{enumerate}

\usepackage[margin=1.2in]{geometry}
\usepackage[english]{babel}
\DeclareMathOperator{\dimH}{dim_H}

\DeclareMathOperator{\Ker}{Ker}

\newcommand{\RR}{\mathbb{R}}
\newcommand{\NN}{\mathbb{N}}
\newcommand{\PP}{\mathbb{P}}
\newcommand{\EE}{\mathbb{E}}

\newenvironment{pr}{\par\noindent {\textit{Proof.}}}
{\begin{flushright} \vspace*{-6mm}\mbox{$\Box$} \end{flushright}}

\newenvironment{proof1}{\par\noindent {\textit{Proof of Theorem \ref{th:LevyMultiple}.}}}
{\begin{flushright} \vspace*{-6mm}\mbox{$\Box$} \end{flushright}}

\begin{document}

\title{On the double points of operator stable L\'evy processes\thanks{Research of T. Luks was partially supported by 
Agence Nationale de la Recherche grant ANR-09-BLAN-0084-01. Research of Y. Xiao was partially supported by 
the NSF Grants  DMS-1307470 and DMS-1309856.}}

\author{Tomasz Luks \and Yimin Xiao}

\institute{T. Luks  \at
              Ecole Centrale de Marseille, I2M\\
38 Rue Fr\'ed\'eric Joliot Curie, 13013 Marseille, France,\\
and Department of Statistics and Probability, Michigan State University\\
\email{tomasz.luks@centrale-marseille.fr}         
           \and
           Y. Xiao \at
            Department of Statistics and Probability, Michigan State University\\
619 Red Cedar Road, C413 Wells Hall,
East Lansing MI 48824-1027 \\
\email{xiaoyimi@stt.msu.edu}
}

\date{Received: / Accepted: }

\maketitle

\begin{abstract}
We determine the Hausdorff dimension of the set of double points for a symmetric operator
stable L\'evy process  $X=\left\{X(t),t\in\RR_+\right\}$ in terms of the eigenvalues of its stability exponent.
\keywords{Multiple points \and Hausdorff dimension \and Operator stable process \and L\'evy process}
\subclass{60J25 \and 60J30 \and 60G51 \and 60G17}
\end{abstract}

\section{Introduction}\label{sec1}

Let $X=\left\{X(t),t\in\RR_+\right\}$ be a stochastic process with values in $\RR^d$, $d\geq 1$,
and let $k\geq2$ be an integer. A point $x\in\RR^d$ is called a {\it k-multiple point} of $X$ if
there exist $k$ distinct times $t_1,\ldots,\, t_k\in\RR_+$ such that
$$
X(t_1)= \ldots =X(t_k)=x.
$$
Denote by $M_k$  the set of $k$-multiple points of $X$.
If $k=2$, then $x$ is also called a {\it double point} of $X$. The existence of multiple points
(or intersections) has been intensely studied for Brownian motion  and more general L\'evy processes
in the literature (see \cite{DEK1,DEK2,DEKT,E,H,He,LRS,MX,S3} and the references therein).

In the present paper we focus on the Hausdorff dimension of the set of double points for
a symmetric operator stable L\'evy process in $\RR^d$, where $d\geq 2$. A L\'evy process
$X=\left\{X(t),t\in\RR_+\right\}$ with values in $\RR^d$ is called {\it operator stable} 
if the distribution $\nu$ of $X(1)$ is full (i.e., not supported on any $(d-1)$-dimensional 
hyperplane) and there exists a linear operator $B$ on $\RR^d$ such that $\displaystyle\nu^t 
=t^B\nu$ for all $t>0$, where $\nu^t$ denotes the $t$-fold convolution power of the infinitely 
divisible law $\nu$ and $t^B\nu(dx)=\nu(t^{-B}dx)$ is the image measure of $\nu$ under the 
linear operator $t^B$. The operator $B$ is called a {\it stability exponent} of $X$. We 
refer to \cite{MS} for more information on operator stable laws.

Our first result gives a general formula for the Hausdorff dimension of the set of 
$k$-multiple points for a symmetric L\'evy process  $X =\left\{X(t),t\in\RR_+\right\}$ 
in terms of its characteristic exponent $\Psi$. See Section 2 for the terminology.

\begin{theorem}\label{th:LevyMultiple}
Let $X =\left\{X(t),t\in\RR_+\right\}$ be a symmetric, absolutely continuous L\'evy process
with values in $\RR^d$ ($d \ge 1$) and L\'evy exponent $\Psi$. Then for any $k\geq2$ we have
$$
\dimH M_k= d-\inf\left\{\beta\in(0,d):\int_{\RR^{kd}}\left[\frac{1}{1+\|\sum_{l=1}^k\xi_l\|^{\beta}}
\prod_{j=1}^k\frac{1}{1+\Psi(\xi_j)} \right]d\overline{\xi}<\infty \right\}
$$
almost surely, where $\overline{\xi}:=(\xi_1,...,\xi_k)\in\RR^{kd}$ and $\xi_i\in\RR^d$ for $i=1,...,k$. 
If the integral above is infinite for all $\beta\in(0,d)$, then $\dimH M_k=0$ a.s.
\end{theorem}

\noindent Next we apply Theorem~\ref{th:LevyMultiple}  to  a symmetric operator
stable L\'evy process with stability exponent $B$. To this end, we factor the minimal
polynomial of $B$ into $q_1(x)\cdots q_p(x)$, where all the roots of $q_i(x)$ have real
parts $a_i$ and $a_i<a_j$ for $i<j$. Define $V_i=\Ker(q_i(B))$ and $d_i=\dim(V_i)$.
Then $d_1+ \cdots +d_p=d$ and $V_1\oplus\cdots\oplus V_p$ is a direct sum decomposition
of $\RR^d$ into $B$-invariant subspaces. We may write $B=B_1\oplus\cdots\oplus B_p$,
where $B_i\colon V_i\to V_i$ and every eigenvalue of $B_i$ has real part equal to $a_i$.
For $j=1,...,d$ and $l=1,...,p$ denote $\alpha_j=a_l^{-1}$ whenever $\sum_{i=0}^{l-1}d_i<j
\leq\sum_{i=0}^{l}d_i$, where $d_0:=0$.  We then have $\alpha_1\geq...\geq\alpha_d$,
and note that $0<\alpha_j\leq2$ in view of \cite[Theorem 7.2.1]{MS}.
Our second theorem, which is the main result of the paper, provides an explicit formula
for the Hausdorff dimension of $M_2$ in terms of the exponents $\alpha_j$.

\begin{theorem}\label{th:OperatorDouble}
Let $X =\left\{X(t),t\in\RR_+\right\}$ be a symmetric operator stable L\'evy process
in $\RR^d$ with stability exponent $B$ and let $M_2$ be the set of double points of $X$. 
\begin{itemize}
\item[(a)]\, If $d=2$ then 
$$
\dimH M_2=\min\left\{\alpha_1\left(2-\frac{1}{\alpha_1}-\frac{1}{\alpha_2}\right), \,
2\alpha_2-\frac{2\alpha_2}{\alpha_1}\right\} \ \ \hbox{a.s.}
$$
\item[(b)]\, If $d=3$ then
$$
\dimH M_2=\alpha_1\left(2-\frac{1}{\alpha_1}-\frac{1}{\alpha_2}-\frac{1}{\alpha_3}\right) \ \ \hbox{a.s.},
$$
\end{itemize}
where a negative dimension means that $M_2=\emptyset$. Furthermore, $M_2=\emptyset$ 
for all $d\geq4$.
\end{theorem}

Theorem~\ref{th:OperatorDouble} is more general than Theorem 1 in \cite{S3}, where $B$ is 
assumed to be a diagonal matrix with entries on the diagonal $\alpha_j \in (1, 2)$ ($1 \le j \le d$). 
Also, the methods used in \cite{S3} are probabilistic in nature and they provide formulas for 
all $k\geq2$. Our approach is analytical and it extends the results of \cite{S3} to the whole 
family of  symmetric operator stable L\'evy processes in the case $k=2$. In addition, we 
provide a necessary and sufficient condition for the existence of double points of symmetric 
operator  stable L\'evy processes in terms of the exponents $\alpha_i$, see Theorem~\ref{th:existd2} 
and Theorem~\ref{th:existd3} below. The latter also  reveals some subtle behavior when $B$ is 
not a diagonal matrix. Adapting our techniques to the case $k\geq3$ is more involved and will 
be dealt with separately.  

The paper is organized as follows. In Section 2 we give definitions and some facts concerning 
operator stable L\'evy processes and we prove Theorem~\ref{th:LevyMultiple}. In Sections 3 and 
4 we focus on the double points problems for $d=2$ and $d=3$, respectively.

Throughout the rest of the paper, $C$ will denote a positive constant, whose value may change 
in each appearance.

\section{Preliminaries}\label{sec2}

A stochastic process $X=\left\{X(t),t\in\RR_+\right\}$ with values in $\RR^d$ is called a 
L\'evy process if $X$ has stationary and independent increments, $X(0)=0$ a.s. and 
$t\mapsto X(t)$ is continuous in probability. We refer to the books \cite{B96,S99} 
for systematic accounts on L\'evy processes.

It is known that the finite-dimensional distributions of $X$ are determined by 
the characteristic function
$$
\EE[ e^{i \langle \xi,\, X(t)\rangle }]= e^{-t\Psi(\xi)}, \quad 
\forall\, t\geq0,
$$
where $\Psi\colon\RR^d\mapsto\mathbb {C}$ is given by the 
L\'evy-Khintchine formula and is called the {\it characteristic or L\'evy exponent} of $X$. 

A L\'evy process $X$ is said to be 
{\it symmetric} if $-X$ and $X$ have the same finite-dimensional distributions. 
In such a case,  $\Psi (\xi) \ge 0$ for all $\xi \in \RR^d$. Using the
terminology in \cite{KX4,KX5}, we say that $X$ is {\it absolutely continuous}, 
if for all $t>0$, the function $\xi\mapsto e^{-t\Psi(\xi)}$ is in $L^1(\RR^d)$.

The following theorem gives a necessary and sufficient condition for the existence 
of $k$-multiple points for a symmetric, absolutely continuous L\'evy process in 
terms of its characteristic exponent. See \cite{E,LRS,FS,KX3} for appropriate 
conditions in terms of the potential density.

\begin{theorem}\label{th:MultipleExistence}
Let $X =\left\{X(t),t\in\RR_+\right\}$ be a symmetric, absolutely continuous 
L\'evy process with values in $\RR^d$ ($d \ge 1$) and characteristic exponent $\Psi$.
Then $X$ has $k$-multiple points if and only if
$$
\int_{\RR^{d(k-1)}}\prod_{j=1}^k\frac{1}{1+\Psi(\xi_{j-1}-\xi_j)} d\overline{\xi}<\infty,
$$
where $\overline{\xi}:=(\xi_1,...,\xi_{k-1})\in\RR^{d(k-1)}$, $\xi_i\in\RR^d$ for $i=1,...,k-1$, and $\xi_0=\xi_k:=0$.
\end{theorem}

\begin{pr}
The existence of $k$-multiple points of $X$ is equivalent with the existence 
of intersections of $k$ independent copies of $X$, see \cite[Proof of Theorem 1]
{LRS}. Furthermore, since $X$ is symmetric, by \cite[Theorem 2.1]{KX5} $X$ 
is also weakly unimodal. Hence, by \cite[Remark 6.6]{KX4}, $k$ independent 
copies of $X$ intersect if and only if
$$
\int_{\RR^{d(k-1)}}\frac{1}{1+\Psi(\sum_{j=1}^{k-1}v_j)}
\prod_{j=1}^{k-1}\frac{1}{1+\Psi(v_j)} d\overline{v}<\infty.
$$
Applying the change of variables $v_j=\xi_j-\xi_{j+1}$, $j=1,...,k-1$, 
where $\xi_k:=0$, we obtain the desired result.
\end{pr}

According to Theorem~\ref{th:MultipleExistence}, a symmetric, absolutely 
continuous L\'evy process $X$ in $\RR^d$ with L\'evy exponent $\Psi$ has 
double points if and only if
$$
\int_{\RR^d}\left(\frac{1}{1+\Psi(\xi)}\right)^2d\xi<\infty.
$$
Since $\Psi(\xi)\leq\|\xi\|^2$ for all $\|\xi\|$ large enough, the above condition 
does not hold when $d\geq4$, which immediately implies the last statement 
of Theorem~\ref{th:OperatorDouble}.\\

\begin{proof1}
Let $X$ be a symmetric, absolutely continuous L\'evy process in $\RR^d$ and 
let $X_1,\, \ldots,\,X_k$ be $k$ independent copies of $X$. For any  $x_1,x_2,\, ...,\,
x_k \in \RR^d$, let $\widetilde{X}_j(t) = x_j + X_j(t)$ for all $t \ge0$ and $1 \le j \le k$. 
Thus, $\widetilde{X}_j = \{\widetilde{X}_j(t), t \ge 0\}$ is a L\'evy process starting 
from $x_j$. 

Denote by $\widetilde{M}_k$ the set of intersections of $\widetilde{X}_1$,\,...,
$\widetilde{X}_k$, i.e.,
$$
\widetilde{M}_k=\bigcap_{i=1}^{k} \widetilde{X}_i\left([0,\infty)\right).
$$
Each $X_i$ has a one-potential density $u\colon\RR^d\mapsto\overline{\RR}_+$ 
satisfying $u(0)>0$. Indeed, by the symmetry, the transition density of $X_i$ 
satisfies
$$
p_t(0)=\int_{\RR^d}\left(p_{t/2}(x)\right)^2dx>0.
$$
By \cite[Theorem 1.5]{KX3}, we can derive that, for $x_1, \ldots, x_k$ in a small 
neighborhood of the origin, the formula of Theorem~\ref{th:LevyMultiple} holds 
for $\dimH\widetilde{M}_k$. Therefore, it is enough to show that $\dimH M_k= 
\dimH\widetilde{M}_k$ almost surely. This last property is part of the folklore 
for L\'evy processes and has already been applied by  
e.g., \cite[p. 85]{H} or \cite[Section 3]{S3}. However, as pointed out in 
\cite[p. 510]{LRS}, the equality needs some rigorous justification in our 
general context.

Since the idea is similar to the argument in the proof of \cite[Theorem 1]{LRS}, 
we only provide the main steps of our proof. 
For a positive integer $N$ and $\alpha \in(0,2)$ let $Y_{\alpha,N}$ be an 
$N$-parameter additive $\alpha$-stable process with values in $\RR^d$ (cf. \cite{KX4,KX3}), 
independent of  $X,\,X_1$,...,$X_k$. 
By \cite[Corollary 3.4, see also Section 4]{KX1}, for any Borel set $F\subset\RR^d$ we have
\begin{equation}\label{eq:Hdimension}
\dimH F=d-\inf\left\{N\alpha>0: \PP(F\cap Y_{\alpha,N}(\RR_+^N))>0  \right\}.
\end{equation}
For $0\leq s<t$,  let $X([s,t])$ denotes the path of $X$ on the time interval $[s,t]$.
The following conditions are equivalent.
\begin{enumerate}
\item[(i)] $\PP(\widetilde{M}_k\cap Y_{\alpha,N}(\RR_+^N))>0$.
\item[(ii)] There exist a constant $s>0$ and a neighborhood $U$ of 0 in $\RR^d$ such 
that for any (initial states of $\widetilde{X}_i$) $x_1,..., x_k \in U$ and $r>0$
we have
$$
\PP\bigg(Y_{\alpha,N}(\RR_+^N)\cap\bigcap_{i=1}^{k} \widetilde{X}_i([0,s])\neq\emptyset, 
\widetilde{X}_i(s)\in B(0,r),\, i=1,...,k \bigg)>0.
$$
\item[(iii)] There exist a constant $s>0$ and positive real numbers $a_1,\,...\,, a_{k-1}$ 
satisfying $a_i\geq a_{i-1}+2$, $a_0=0$, such that
$$
\PP\bigg(Y_{\alpha,N}(\RR_+^N)\cap\bigcap_{i=0}^{k-1} X([a_is,\, (a_i+1)s])\neq\emptyset\bigg)>0.
$$
\end{enumerate}
The equivalence between (i) and (ii) follows from the fact that $p_t(0)>0$ for all $t>0$. Indeed, by 
the lower semicontinuity of $p_t$, for any $t>0$ there exists a neighborhood $U$ of 0 such that for 
all $x,y\in U$ we have $p_t(x-y)>0$, see also \cite[Proof of Theorem 3]{LRS}. The equivalence 
between (ii) and (iii) can be proved using similar techniques as in \cite[Proof of Theorem 1]{LRS}. 
Finally, Condition (iii) is equivalent with $\PP(M_k\cap Y_{\alpha,N}(\RR_+^N))>0$. This,  
(\ref{eq:Hdimension}) and the Borel-Cantelli argument in  \cite[p. 511]{LRS} yield Theorem 
\ref{th:LevyMultiple}.\\
\end{proof1}

Let $X$ be an operator stable L\'evy process in $\RR^d$ with stability exponent $B$, where $d \ge 2$.
It is well-known that there exists a real invertible $d\times d$ matrix $P$ such that
$$
B=PDP^{-1},
$$
where $D$ is a real $d\times d$ matrix of the form
\begin{equation}\label{eq:matrix1}
	D=\left(
	\begin{matrix}
		J_1 & 0 &  \ldots & 0 \\
		0 & J_2 &  \ldots & 0 \\
		\vdots & \vdots  & \ddots & \vdots \\
		0 & 0 & \ldots  & J_p 	
		\end{matrix}
	\right),
\end{equation}
and each block $J_i$, $i=1,...,p$, is of the form

\begin{equation}\label{eq:matrix2}
	\left(
	\begin{matrix}
		a & 0 & 0 & \ldots & 0 \\
		1 & a & 0 & \ldots & 0 \\
		0 & 1 & a & \ldots & 0 \\
		\vdots &   & \ddots & \ddots & \vdots \\
		0 & \ldots  & 0 & 1 & a
	\end{matrix}
	\right)
	\quad \hbox{ or } \quad
	\left(
	\begin{matrix}
		C & 0 & 0 & \ldots & 0 \\
		I & C & 0 & \ldots & 0 \\
		0 & I & C & \ldots & 0 \\
		\vdots &   & \ddots & \ddots & \vdots \\
		0 & \ldots  & 0 & I & C
	\end{matrix}
	\right),
\end{equation}
where $a$ is a real eigenvalue of $B$ in the first case, and in the second case
\begin{equation}\label{eq:matrix3}
	C=\left(
	\begin{matrix}
		a & -b \\
		b & a  \\
	\end{matrix}
	\right)
	\quad{\rm  and} \quad
	I=\left(
	\begin{matrix}
		1 & 0 \\
		0 & 1  \\
	\end{matrix}
	\right),
\end{equation}
where $a\pm ib$ is a complex conjugate pair of eigenvalues of $B$. Using
the notation from the Introduction, the size of $J_i$ is equal to $d_i$ and
$a_i$ is the real part of the corresponding eigenvalue of $J_i$. Recall
also that $\alpha_j=a_l^{-1}$ whenever $\sum_{i=0}^{l-1}d_i<j
\leq\sum_{i=0}^{l}d_i$, where $d_0:=0$,  $j=1,...,d$ and $l=1,...,p$.

If $X$ is symmetric, then the L\'evy exponent $\Psi$ of $X$ is nonnegative,
and by \cite[Theorem 4.2]{MX} the following condition holds: for every
$\varepsilon>0$, there exists a constant $\tau>1$ such that
\begin{equation}\label{est:Levyexp}
\frac{K^{-1}}{\|\xi\|^{\varepsilon}\sum_{j=1}^d|\xi_j|^{\alpha_j}}\leq
\frac{1}{1+\Psi(\xi)}\leq\frac{K\|\xi\|^{\varepsilon}}{\sum_{j=1}^d|\xi_j|^{\alpha_j}}
\end{equation}
for all $\xi\in\RR^d$ with $\|\xi\|\geq\tau$, where $K\geq1$ is a
constant which depends on $\varepsilon$ and $\tau$ only.

\section{Double points problem for $d=2$}\label{sec3}

Throughout this section, $X$ is a symmetric operator stable L\'evy process
in $\RR^2$ with L\'evy exponent $\Psi$ and stability exponent $B$ whose eigenvalues
have real parts $\alpha_1^{-1},\alpha_2^{-1}$, as explained in Section~\ref{sec1}.

To prove Part (a) of Theorem~\ref{th:OperatorDouble},
we apply Theorem~\ref{th:LevyMultiple} to $X$  using the estimate (\ref{est:Levyexp}).
Let
$$
A=\left\{(x,y)\in\RR^2\times\RR^2:|x|>1, |y|>1\right\},
$$
and let
$$
I_{\beta}=\iint_{A}\frac{dxdy}{(1+|x_1+y_1|^{\beta}+|x_2+y_2|^{\beta})
(|x_1|^{\alpha_1}+|x_2|^{\alpha_2})(|y_1|^{\alpha_1}+|y_2|^{\alpha_2})}.
$$
By Theorem~\ref{th:LevyMultiple}, (\ref{est:Levyexp}) and by Lemma~\ref{lem:Ibeta} below, we have
\begin{equation}\label{eq:Hausdorff1}
\dimH M_2= 2-\inf\left\{\beta\in(0,2):I_{\beta}<\infty \right\}  \ \ \hbox{a.s.}
\end{equation}
Therefore, Part (a) of Theorem~\ref{th:OperatorDouble} is an immediate
consequence of (\ref{eq:Hausdorff1}) and Theorem~\ref{th:Ibeta} below.

\begin{theorem}\label{th:Ibeta}
If $2-1/\alpha_1-1/\alpha_2>0$ then
\begin{equation} \label{Eq:Th4}
\inf\left\{\beta>0:I_{\beta}<\infty\right\}=\max\left\{3+\frac{\alpha_1}{\alpha_2}-
2\alpha_1,\, 2+\frac{2\alpha_2}{\alpha_1}-2\alpha_2\right\}.
\end{equation}
If $2-1/\alpha_1-1/\alpha_2\leq0$, then $I_{\beta}=\infty$ for any $\beta>0$.
Thus $\dimH M_2=0$ a.s.
\end{theorem}

\begin{pr}
We divide the proof into two parts. Part (i) considers the case when both $\alpha_1$ 
and $\alpha_2$ are not integers. Part (ii) deals with the remaining cases [i.e., $\alpha_1 =  \alpha_2  = 2$ 
and $\alpha_1 = 2, \alpha_2 = 1$], where an extra factor of $\ln k$ or 
$\ln n$ may appear in upper bounds for the integrals in (\ref{eq:int1})--(\ref{eq:int4}) 
below. Since only slight modifications will be needed to prove (\ref{Eq:Th4}),  we will not provide
all the details.

{\it Part (i):}  For $k,n\in\NN$ denote by $A_{k,n}$ the set
\[
\begin{split}
\Big\{(x,y)\in\RR^2\times\RR^2:\,|x|>1, |y|>1,\, k-1\leq|x_1+y_1|<k,\, n-1\leq|x_2+y_2|<n\Big\}.
\end{split}
\]
We have $I_{\beta}<\infty$ if and only if
\begin{equation}\label{eq:2series}
\sum_{k=1}^{\infty}\sum_{n=1}^{\infty}\frac{1}{k^{\beta}+n^{\beta}}
\int\int_{A_{k,n}}\frac{dxdy}{(|x_1|^{\alpha_1}+|x_2|^{\alpha_2})(|y_1|^{\alpha_1}+|y_2|^{\alpha_2})}
<\infty.
\end{equation}
Furthermore, we can assume that $|x_1|, |x_2| \geq1$ and $y_1,y_2\geq1$. We then have
$$
|y_1-k|\wedge|y_1-k+1|\leq|x_1|\leq y_1+k,
$$
$$
|y_2-n|\wedge|y_2-n+1|\leq|x_2|\leq y_2+n.
$$
The first estimate follows from the inequalities
$$
-y_1-k\leq x_1\leq-y_1-k+1,
$$
or
$$
-y_1+k-1\leq x_1\leq-y_1+k,
$$
for the second one the argument is similar. Furthermore, we may assume 
that $|x_1|\geq |y_1-k|$, $|x_2|\geq |y_2-n|$ and $|y_1-k|+|y_2-n|\geq C>0$.
Therefore, the double series in (\ref{eq:2series}) is bounded from below by
\begin{equation}\label{eq:2seriesLow}
\sum_{k=1}^{\infty}\sum_{n=1}^{\infty}\frac{1}{k^{\beta}+n^{\beta}}
\int_1^{\infty}\int_1^{\infty}\frac{dy_1dy_2}{((y_1+k)^{\alpha_1}+(y_2+n)^{\alpha_2})
(y_1^{\alpha_1}+y_2^{\alpha_2})},
\end{equation}
and from above by
\begin{equation}\label{eq:2seriesUp}
\sum_{k=1}^{\infty}\sum_{n=1}^{\infty}\frac{1}{k^{\beta}+n^{\beta}}
\int_{B_{k,n}}\frac{dy}{(|y_1-k|^{\alpha_1}+|y_2-n|^{\alpha_2})(y_1^{\alpha_1}+y_2^{\alpha_2})},
\end{equation}
where $$B_{k,n}=\left\{y\in\RR^2:\min\left\{y_1,y_2,|y_1-k|,|y_2-n|\right\}\geq1\right\}.$$
To estimate (\ref{eq:2seriesUp}), we may assume that $k,n\geq3$ and note that
\begin{align}
&\int_{B_{k,n}}\frac{dy}{(|y_1-k|^{\alpha_1}+|y_2-n|^{\alpha_2})(y_1^{\alpha_1}+y_2^{\alpha_2})}\nonumber\\
&=\int_1^{n-1}\int_1^{k-1}\frac{dy_1dy_2}{((k-y_1)^{\alpha_1}+(n-y_2)^{\alpha_2})(y_1^{\alpha_1}+y_2^{\alpha_2})}\label{eq:int1}\\
&+\int_{n+1}^{\infty}\int_1^{k-1}\frac{dy_1dy_2}{((k-y_1)^{\alpha_1}+(y_2-n)^{\alpha_2})(y_1^{\alpha_1}+y_2^{\alpha_2})}\label{eq:int2}\\
&+\int_1^{n-1}\int_{k+1}^{\infty}\frac{dy_1dy_2}{((y_1-k)^{\alpha_1}+(n-y_2)^{\alpha_2})(y_1^{\alpha_1}+y_2^{\alpha_2})}\label{eq:int3}\\
&+\int_{n+1}^{\infty}\int_{k+1}^{\infty}\frac{dy_1dy_2}{((y_1-k)^{\alpha_1}+(y_2-n)^{\alpha_2})(y_1^{\alpha_1}+y_2^{\alpha_2})}.\label{eq:int4}
\end{align}
We start with the integral (\ref{eq:int4}). After a change of variables, (\ref{eq:int4}) is equal to
$$
\int_{1}^{\infty}\int_{1}^{\infty}\frac{dtds}{(t^{\alpha_1}+s^{\alpha_2})((t+k)^{\alpha_1}+(s+n)^{\alpha_2})}.
$$
Observe that the same double integral appears in the lower estimate (\ref{eq:2seriesLow}). We have
$$
\frac{1}{(t+k)^{\alpha_1}+(s+n)^{\alpha_2}}\asymp\frac{1}{t^{\alpha_1}+k^{\alpha_1}+s^{\alpha_2}+n^{\alpha_2}}
\asymp \frac{1}{k^{\alpha_1}+n^{\alpha_2}} \wedge \frac{1}{t^{\alpha_1}}\wedge\frac{1}{s^{\alpha_2}}.
$$
Here $\asymp$ means that there are two sided estimates with a constant depending only on $\alpha_1$ and $\alpha_2$. 
Hence, (\ref{eq:int4}) is comparable with
\begin{align*}
&\int_{1}^{\infty}\int_{1}^{\infty}\left(\frac{1}{t^{\alpha_1}}\wedge\frac{1}{s^{\alpha_2}}\right)
\left(\frac{1}{k^{\alpha_1}+n^{\alpha_2}} \wedge \frac{1}{t^{\alpha_1}}
\wedge\frac{1}{s^{\alpha_2}} \right)\,dtds\\
&=\int_{1}^{\infty}\int_{1}^{s^{\alpha_2/\alpha_1}}\frac{1}{s^{\alpha_2}}
\left(\frac{1}{k^{\alpha_1}+n^{\alpha_2}} \wedge\frac{1}{s^{\alpha_2}}\right) \,dtds\\
&+\int_{1}^{\infty}\int_{s^{\alpha_2/\alpha_1}}^{\infty}\frac{1}{t^{\alpha_1}}
\left(\frac{1}{k^{\alpha_1}+n^{\alpha_2}} \wedge\frac{1}{t^{\alpha_1}}\right)\,dtds 
=:I_1+I_2.
\end{align*}
We write
\begin{align*}
I_2=&\int_{1}^{\infty}\int_{s^{\alpha_2/\alpha_1}}^{\infty}\frac{1}{t^{\alpha_1}}
\left(\frac{1}{k^{\alpha_1}+n^{\alpha_2}} \wedge\frac{1}{t^{\alpha_1}}\right)\, dtds\\
=&\int_{1}^{(k^{\alpha_1}+n^{\alpha_2})^{1/\alpha_2}}\int_{s^{\alpha_2/\alpha_1}}^{\infty}
\frac{1}{t^{\alpha_1}} \left(\frac{1}{k^{\alpha_1}+n^{\alpha_2}} \wedge\frac{1}{t^{\alpha_1}}\right)\,dtds\\
&+\int^{\infty}_{(k^{\alpha_1}+n^{\alpha_2})^{1/\alpha_2}}\int_{s^{\alpha_2/\alpha_1}}^{\infty}
\frac{1}{t^{2\alpha_1}}dtds =:I_2^{(1)}+I_2^{(2)}.
\end{align*}
It can be seen that if $2-1/\alpha_1-1/\alpha_2>0$, then 
$$
I_2^{(2)}=C\int^{\infty}_{(k^{\alpha_1}+n^{\alpha_2})^{1/\alpha_2}}s^{\alpha_2/\alpha_1-2\alpha_2}ds=
C\left(k^{\alpha_1}+ n^{\alpha_2}\right)^{1/\alpha_1+1/\alpha_2-2}.
$$
On the other hand, if $2-1/\alpha_1-1/\alpha_2\leq0$, then
the integral (\ref{eq:int4}) is infinite and so is the series (\ref{eq:2seriesLow}). 
So  $I_{\beta}=\infty$ for any $\beta>0$, which proves the second part of the theorem.

Next we consider the case when $2-1/\alpha_1-1/\alpha_2>0$. This implies that $\alpha_1>1$. 
We have
\begin{align*}
I_2^{(1)}&=\frac{1}{k^{\alpha_1}+n^{\alpha_2}}\int_{1}^{(k^{\alpha_1}+n^{\alpha_2})^{1/\alpha_2}}\int_{s^{\alpha_2/\alpha_1}}^{(k^{\alpha_1}+n^{\alpha_2})^{1/\alpha_1}}t^{-\alpha_1}dtds\\
&\qquad \quad +\int_{1}^{(k^{\alpha_1}+n^{\alpha_2})^{1/\alpha_2}}ds\int^{\infty}_{(k^{\alpha_1}+n^{\alpha_2})^{1/\alpha_1}}t^{-2\alpha_1}dt\\
&\leq\frac{C}{k^{\alpha_1}+n^{\alpha_2}}\int_{1}^{(k^{\alpha_1}+n^{\alpha_2})^{1/\alpha_2}}s^{\alpha_2/\alpha_1-\alpha_2}ds
+C(k^{\alpha_1}+n^{\alpha_2})^{1/\alpha_1+1/\alpha_2-2}\\
&\leq C(k^{\alpha_1}+n^{\alpha_2})^{1/\alpha_1+1/\alpha_2-2}.
\end{align*}
In the last inequality we have used the assumption $\alpha_2<2$. Similarly, we verify that
\begin{align*}
I_1&\leq\int_{1}^{\infty}s^{\alpha_2/\alpha_1-\alpha_2} \left(\frac{1}{k^{\alpha_1}
+n^{\alpha_2}} \wedge\frac{1}{s^{\alpha_2}}\right)\, ds\\
&=\int_{1}^{(k^{\alpha_1}+n^{\alpha_2})^{1/\alpha_2}}\frac{s^{\alpha_2/\alpha_1-\alpha_2}}{k^{\alpha_1}+n^{\alpha_2}}ds
+\int_{(k^{\alpha_1}+n^{\alpha_2})^{1/\alpha_2}}^{\infty}s^{\alpha_2/\alpha_1-2\alpha_2}ds\\
&\leq C\left(k^{\alpha_1}+ n^{\alpha_2}\right)^{1/\alpha_1+1/\alpha_2-2}.
\end{align*}
Hence we have
\begin{equation}\label{Eq:14}
(\ref{eq:int4})\asymp\left( \frac 1 {k^{\alpha_1}+ n^{\alpha_2}}\right)^{2- 1/\alpha_1- 1/\alpha_2}.
\end{equation}

Now we consider the integral (\ref{eq:int1}). After the change of variables $t=k-y_1$, $s=n-y_2$, we get
\[
\begin{split}
&\int_1^{\frac{n}{2}}\int_1^{k-1}\frac{dy_1dy_2}{\big[(k-y_1)^{\alpha_1}+(n-y_2)^{\alpha_2}\big](y_1^{\alpha_1}+y_2^{\alpha_2})}\\
&=\int_{\frac{n}{2}}^{n-1}\int_1^{k-1}\frac{dtds} {(t^{\alpha_1}+s^{\alpha_2}) \big[(k-t)^{\alpha_1}+(n-s)^{\alpha_2}\big]}.
\end{split}
\]
Hence, the integral (\ref{eq:int1}) is equal to
\begin{equation}
\begin{split}\label{Eq:10}
&2\int_1^{\frac{n}{2}}\int_1^{k-1}\frac{dtds}{(t^{\alpha_1}+s^{\alpha_2})\big[(k-t)^{\alpha_1}+(n-s)^{\alpha_2}\big]}
\leq \frac{C}{n^{\alpha_2}}\int_1^{\frac{n}{2}}\int_1^{k-1}\frac{dtds}{(t+s^{\alpha_2/\alpha_1})^{\alpha_1}}\\
&=\frac{C}{n^{\alpha_2}}\int_1^{\frac{n}{2}}\int_{1+s^{\alpha_2/\alpha_1}}^{k+s^{\alpha_2/\alpha_1}-1}u^{-\alpha_1}duds
\leq \frac{C}{n^{\alpha_2}}\int_1^{\frac{n}{2}}s^{\alpha_2/\alpha_1-\alpha_2}ds,
\end{split}
\end{equation}
where the last inequality follows from the fact that $\alpha_1>1$. Since $\alpha_2<2$, we get
\begin{equation}\label{alpha}
\frac{C}{n^{\alpha_2}}\int_1^{\frac{n}{2}}s^{\alpha_2/\alpha_1-\alpha_2}ds
\leq C\left( \frac 1 {n^{\alpha_2}}\right)^{2- 1/\alpha_1 -1/\alpha_2}.
\end{equation}
Interchanging the roles of $n$ and $k$, the same argument as in (\ref{Eq:10}) 
and (\ref{alpha}) shows that the integral (\ref{eq:int1}) is at most
$$
C\left(\frac{1}{k^{\alpha_1}}\wedge\frac{1}{n^{\alpha_2}}\right)^{2-1/\alpha_1-1/\alpha_2}
\leq C\left(\frac{1}{k^{\alpha_1}+ n^{\alpha_2}}\right)^{2-1/\alpha_1-1/\alpha_2}.
$$
For the integral (\ref{eq:int2}), a change of variables $k-y_1=t$, $y_2-n=s$ implies that 
it  is equal to
\begin{align*}
&\int_1^{\infty}\int_1^{k-1}\frac{dtds}{(t^{\alpha_1}+s^{\alpha_2})\big[(k-t)^{\alpha_1}+(s+n)^{\alpha_2}\big]}\\
&\leq C\, \int_1^{\infty}\frac{ds}{(s+n)^{\alpha_2}}\int_1^{k-1}\frac{dt}{(t+s^{\alpha_2/\alpha_1})^{\alpha_1}}\\
&\leq C\,\int_1^{\infty}\frac{ds}{(s+n)^{\alpha_2}s^{\alpha_2-\alpha_2/\alpha_1}}
\le C\left(\frac 1 {n^{\alpha_2}}\right)^{2 - 1/\alpha_1 -1/\alpha_2}.
\end{align*}
In deriving the last inequality we have used the assumption that $\alpha_2<2$. 

To bound the integral (\ref{eq:int2}) in terms of $k$, we note that 
$$
\int_1^{\infty}\int_1^{\frac{k}{2}}\frac{dtds}{(t^{\alpha_1}+s^{\alpha_2})\big[(k-t)^{\alpha_1}+(s+n)^{\alpha_2}\big]}
\leq \int_1^{\infty}\int_1^{\frac{k}{2}}\frac{Cdtds}{(t^{\alpha_1}+s^{\alpha_2})(k^{\alpha_1}+s^{\alpha_2})}
$$
and
\begin{align*}
\int_1^{\infty}\int_{\frac{k}{2}}^{k-1}\frac{dtds}{(t^{\alpha_1}+s^{\alpha_2})\big[(k-t)^{\alpha_1}+(s+n)^{\alpha_2}\big]}
\leq \int_1^{\infty}\int_1^{\frac{k}{2}}\frac{Cduds}{(k^{\alpha_1}+s^{\alpha_2})(u^{\alpha_1}+s^{\alpha_2})}.
\end{align*}
Thus,
$$
\int_1^{\infty}\int_1^{k-1}\frac{dtds}{(t^{\alpha_1}+s^{\alpha_2})\big[(k-t)^{\alpha_1}+(s+n)^{\alpha_2}\big]}
\leq \int_1^{\frac{k}{2}}\int_1^{\infty}\frac{Cdsdt}{(t^{\alpha_1}+s^{\alpha_2})(k^{\alpha_1}+s^{\alpha_2})}.
$$
Assume first that $\alpha_2>1$. We can verify that the integral (\ref{eq:int2}) is less than
\begin{align*}
\frac{C}{k^{\alpha_1}}\int_1^{\frac{k}{2}}\int_1^{\infty}\frac{dsdt}{(t^{\alpha_1/\alpha_2}+s)^{\alpha_2}}
\leq C \left(\frac 1 {k^{\alpha_1}}\right)^{2- 1/\alpha_1-1/\alpha_2},
\end{align*}
Combining the above yields that (\ref{eq:int2}) is bounded from above by
$$C\left(k^{\alpha_1}+ n^{\alpha_2}\right)^{1/\alpha_1+1/\alpha_2-2}.$$ 

By symmetry we also get that (\ref{eq:int3}) is less than $C\left(k^{\alpha_1}+ n^{\alpha_2}\right)^{1/\alpha_1+1/\alpha_2-2}$
 (indeed, this case is even easier since $\alpha_1>1$). Therefore, we have proved that
$$
(\ref{eq:int1}),\, (\ref{eq:int2}),\, (\ref{eq:int3})\leq C \left(\frac 1 {k^{\alpha_1}+ n^{\alpha_2}}\right)^{2 - 1/\alpha_1- 1/\alpha_2}.
$$
This and (\ref{Eq:14}) imply that, except the cases $\alpha_1=\alpha_2=2$ or $\alpha_1=2$ and $\alpha_2=1$,  
the series (\ref{eq:2seriesLow})  and (\ref{eq:2seriesUp}) (and hence $I_{\beta}$) 
are finite if and only if $2-1/\alpha_1-1/\alpha_2>0$ and
$$
\sum_{k=1}^{\infty}\sum_{n=1}^{\infty}\frac{1}{k^{\beta}+n^{\beta}}\left(\frac{1}{k^{\alpha_1}+ 
n^{\alpha_2}}\right)^{2-1/\alpha_1-1/\alpha_2}<\infty.
$$
Hence the theorem follows from Lemma~\ref{lem:series} below.

{\it Part (ii):} If $\alpha_2=\alpha_1=2$ then the methods used in Part (i) still apply. 
In this case, the left hand side of (\ref{alpha}) is less than  $Cn^{-2}\ln n$, and by 
symmetry, the integral (\ref{eq:int1}) can be estimated by
\begin{align*}
\frac{\ln k }{k^2}\wedge\frac{\ln n}{n^2} 
\leq C\left(\frac{1}{k^{2-\varepsilon}+n^{2-\varepsilon}}\right)^{2-2/(2-\varepsilon)},
\end{align*}
for any small $\varepsilon>0$ and $k,n\geq N_{\varepsilon}>0$, provided $N_{\varepsilon}$ 
is sufficiently large. Therefore, Lemma~\ref{lem:series} can be applied with $\alpha_1
=\alpha_2=2-\varepsilon$ for the upper estimate and with  $\alpha_1=\alpha_2=2$ for 
the lower bound. Since $\varepsilon >0$ is arbitrary, we obtain the desired result. 
In the case $\alpha_1=2$ and $\alpha_2=1$, the reasoning is similar and is omitted.
\end{pr}

\begin{lemma}\label{lem:series}
We have
\[
\begin{split}
&\inf\left\{\beta\in(0,2):\sum_{k=1}^{\infty}\sum_{n=1}^{\infty}\frac{1}{k^{\beta}+n^{\beta}}
\left(\frac{1}{k^{\alpha_1}+ n^{\alpha_2}}\right)^{2-1/\alpha_1-1/\alpha_2}<\infty\right\}\\
&=\max\left\{3+\frac{\alpha_1}{\alpha_2}-2\alpha_1,2+\frac{2\alpha_2}{\alpha_1}-2\alpha_2\right\}.
\end{split}
\]
\end{lemma}

\begin{pr}
The convergence of the series is equivalent with the convergence of the integral
\begin{align*}
&\int_1^{\infty}\int_1^{\infty}\left(\frac{1}{x}\wedge\frac{1}{y}\right)^{\beta}\left(\frac{1}{x^{\alpha_1}}
\wedge\frac{1}{y^{\alpha_2}}\right)^{2-1/\alpha_1-1/\alpha_2}dxdy\\
&=\int_1^{\infty}\int_1^{y^{\alpha_2/\alpha_1}}\frac{1}{y^{\beta}}\left(\frac{1}{y^{\alpha_2}}\right)^{2-1/\alpha_1-1/\alpha_2}dxdy\\
&\qquad +\int_1^{\infty}\int^y_{y^{\alpha_2/\alpha_1}}\frac{1}{y^{\beta}}\left(\frac{1}{x^{\alpha_1}}\right)^{2-1/\alpha_1-1/\alpha_2}dxdy\\
&\qquad +\int_1^{\infty}\int_y^{\infty}\frac{1}{x^{\beta}}\left(\frac{1}{x^{\alpha_1}}\right)^{2-1/\alpha_1-1/\alpha_2}dxdy
=:I_1+I_2+I_3.
\end{align*}
It can be seen that $I_1<\infty$ if and only if
$$
\int_1^{\infty}y^{-\beta-2\alpha_2+2\alpha_2/\alpha_1+1}dy<\infty,
$$
and the last condition is equivalent with $\beta>2+2\alpha_2/\alpha_1-2\alpha_2$. 

Next we consider $I_2$.
If $\alpha_1=\alpha_2$, then $I_2=0$, so assume that $\alpha_1\neq\alpha_2$. If $-2\alpha_1+\alpha_1/\alpha_2=-2$,  then
$$
I_2= (1-\alpha_2/\alpha_1)\int_1^{\infty}y^{-\beta}\ln y \,dy.
$$
So $I_2<\infty$ if and only if $\beta>1=3+\alpha_1/\alpha_2-2\alpha_1$. Suppose now $-2\alpha_1+\alpha_1/\alpha_2\neq-2$. Then we have
$$
I_2=\int_1^{\infty}y^{-\beta}\left(\frac{y^{-2\alpha_1+\alpha_1/\alpha_2+2}-(y^{\alpha_2/\alpha_1})^{-2\alpha_1+\alpha_1/\alpha_2+2}}{-2\alpha_1+\alpha_1/\alpha_2+2}\right)dxdy.
$$
We consider two cases:
\begin{enumerate}
\item[(a)] If $-2\alpha_1+\alpha_1/\alpha_2+2>0$, then
$$
I_2\leq \int_1^{\infty}\frac{y^{-\beta-2\alpha_1+\alpha_1/\alpha_2+2}}{-2\alpha_1+\alpha_1/\alpha_2+2}dxdy,
$$
and the last integral is finite if $\beta>3+\alpha_1/\alpha_2-2\alpha_1$.
\item[(b)] If $-2\alpha_1+\alpha_1/\alpha_2+2<0$, then
$$
I_2\leq \int_1^{\infty}y^{-\beta}\left(\frac{(y^{\alpha_2/\alpha_1})^{-2\alpha_1+\alpha_1/\alpha_2+2}}{-2\alpha_1+\alpha_1/\alpha_2+2}\right)dxdy
=\int_1^{\infty}\frac{y^{-\beta+2\alpha_2/\alpha_1-2\alpha_2+1}}{-2\alpha_1+\alpha_1/\alpha_2+2}dy,
$$
which  is finite if $\beta>2+2\alpha_2/\alpha_1-2\alpha_2$.
\end{enumerate}
Therefore, the condition $\beta>\max\left\{3+\alpha_1/\alpha_2-2\alpha_1,2+\alpha_2/\alpha_1-2\alpha_2\right\}$ implies $I_2<\infty$.

Finally, we consider $I_3$.
A necessary condition for $I_3<\infty$ is $-\beta-2\alpha_1+\alpha_1/\alpha_2+1<-1$. Assuming this we get
$$
I_3=\int_1^{\infty}\frac{y^{-\beta-2\alpha_1+\alpha_1/\alpha_2+2}}{\beta+2\alpha_1-\alpha_1/\alpha_2-2}dy.
$$
Thus $I_3<\infty$ if and only if $\beta>3+\alpha_1/\alpha_2-2\alpha_1$.

Therefore, we have proved that
$$
\beta>\max\left\{3+\alpha_1/\alpha_2-2\alpha_1,2+\alpha_2/\alpha_1-2\alpha_2\right\} \ \Rightarrow\  I_1,I_2,I_3<\infty,
$$
and
$$
I_1,I_3<\infty\  \Rightarrow \ \beta>\max\left\{3+\alpha_1/\alpha_2-2\alpha_1,2+\alpha_2/\alpha_1-2\alpha_2\right\}.
$$
This yields the conclusion of the lemma.
\end{pr}

\begin{lemma}\label{lem:Ibeta}
We have
\begin{equation}\label{Eq:lemma7}
\begin{split}
&\inf\left\{\beta\in(0,2):I_{\beta}<\infty \right\}\\
&=\inf\left\{\beta\in(0,2):\int_{\RR^2}\int_{\RR^2}\frac{dxdy}{(1+\|x+y\|^{\beta})(1+\Psi(x))(1+\Psi(y))}<\infty\right\}.
\end{split}
\end{equation}
\end{lemma}

\begin{pr}
Denote  the first and the second term  in (\ref{Eq:lemma7}) by $\gamma$ and $\gamma'$, respectively.
For any fixed $\beta>\gamma$, we show that if $\varepsilon>0$ is small enough, then
\begin{equation}\label{Eq:18}
\iint_A\frac{\|x\|^{\varepsilon}\|y\|^{\varepsilon}}{(1+\|x+y\|^{\beta})(|x_1|^{\alpha_1}+|x_2|^{\alpha_2})(|y_1|^{\alpha_1}+|y_2|^{\alpha_2})}dxdy<\infty,
\end{equation}
where $A=\left\{(x,y)\in\RR^2\times\RR^2:|x|>1, |y|>1\right\}$. By the upper bound in (\ref{est:Levyexp}), this implies $\gamma'\leq\gamma$. 

To prove  (\ref{Eq:18}), it is enough to show that
$$
\iint_A\frac{(|x_1|^{\varepsilon\alpha_1}+|x_2|^{\varepsilon\alpha_2})(|y_1|^{\varepsilon\alpha_1}+|y_2|^{\varepsilon\alpha_2})}{(1+|x_1+y_1|^{\beta}+|x_2+y_2|^{\beta})(|x_1|^{\alpha_1}+|x_2|^{\alpha_2})(|y_1|^{\alpha_1}+|y_2|^{\alpha_2})}dxdy<\infty
$$
for sufficiently small $\varepsilon>0$. Furthermore, the above integral is comparable to 
$$
 \iint_A\frac{dxdy}{(1+|x_1+y_1|^{\beta}+|x_2+y_2|^{\beta})(|x_1|^{\alpha'_1}+|x_2|^{\alpha'_2})(|y_1|^{\alpha'_1}+|y_2|^{\alpha'_2})}=:I'_{\beta},
$$
where $\alpha'_i:=(1-\varepsilon)\alpha_i$. By Theorem~\ref{th:Ibeta},
$$
\gamma=\max\left\{3+\frac{\alpha_1}{\alpha_2}-2\alpha_1,2+\frac{2\alpha_2}{\alpha_1}-2\alpha_2\right\},
$$
and since $\beta>\gamma$, we may choose $\varepsilon>0$ such that
$$
\max\left\{3+\frac{\alpha'_1}{\alpha'_2}-2\alpha'_1,2+\frac{2\alpha'_2}{\alpha'_1}-2\alpha'_2\right\}<\beta.
$$
Theorem~\ref{th:Ibeta} implies  $I'_{\beta} < \infty$ and thus (\ref{Eq:18}) holds. In order to show that 
$\gamma'\geq\gamma$
we use similar arguments and the lower estimate of (\ref{est:Levyexp}).
\end{pr}

According to Part (a) of Theorem~\ref{th:OperatorDouble}, the set of double points of $X$ has 
positive Hausdorff dimension if and only if  $2- 1/{\alpha_1} -  1/{\alpha_2}>0$. The next theorem 
shows that this is also a necessary condition for the existence of double points of $X$.

\begin{theorem}\label{th:existd2}
$M_2$ is nonempty if and only if $2-1/\alpha_1-1/\alpha_2>0$.
\end{theorem}

\begin{pr}
By Theorem \ref{th:MultipleExistence} we have $M_2\neq\emptyset$ if and only if
\begin{equation}\label{eq:intPsi}
\int_{\RR^2}\bigg(\frac{1}{1+\Psi(\xi)}\bigg)^2d\xi<\infty.
\end{equation}
According to the decomposition described in Section~\ref{sec2}, the stability exponent of 
$X$ satisfies $B=PDP^{-1}$. Since we consider the case $d=2$, the matrix $D$ can have the following forms
\begin{enumerate}
\item[(a)] $\left(
	\begin{matrix}
		a_1 & 0 \\
		0 & a_2  \\
	\end{matrix}
	\right); $
\item[(b)] $\left(
	\begin{matrix}
		a & 0 \\
		1 & a  \\
	\end{matrix}
	\right) $ or $\left(
	\begin{matrix}
		a & -b \\
		b & a  \\
	\end{matrix}
	\right)$.
\end{enumerate}
By \cite[(4.9),(4.14),(4.15),(4.16)]{MX}, we have the following estimates of the L\'evy 
exponent $\Psi(\xi)$ when $\|\xi\|\to\infty$, depending on the cases (a) and (b):
\begin{enumerate}
\item[(a)] $\alpha_1=1/a_1$, $\alpha_2=1/a_2$, and $$\displaystyle\Psi(\xi)\asymp |\xi_1|^{\alpha_1}+|\xi_2|^{\alpha_2};$$
\item[(b)] $\alpha_1=\alpha_2=1/a$, and $$\displaystyle\Psi(\xi)\asymp |\xi_1|^{\alpha_1}+|\xi_2|^{\alpha_1}(\ln\|\xi\|)^{\alpha_1}.$$
\end{enumerate}
In the case (a), it follows from the proof of Theorem~\ref{th:Ibeta} and the estimates of the 
integral (\ref{eq:int4}) that (\ref{eq:intPsi}) holds if and only if $2-1/\alpha_1-1/\alpha_2>0$.
In the case (b) we have $\alpha_1=\alpha_2$ and the inequality $2-1/\alpha_1-1/\alpha_2>0$ 
is equivalent with $\alpha_1>1$. Hence, it is enough to show that (\ref{eq:intPsi}) does not hold
for $\alpha_1=1$. Under this assumption we have
\begin{align*}
\int_{\RR^2}\left(\frac{1}{1+\Psi(\xi)}\right)^2d\xi&\geq C\int_1^{\infty}\int_1^{\infty}\frac{dxdy}{\left(x+y\ln\sqrt{x^2+y^2}\right)^2}\\
&\geq 
C\int_1^{\infty}\int_1^x\frac{dydx}{\left(x+y\ln x\right)^2} \\
&=\int_1^{\infty}\left(\frac{1}{x+\ln x}-\frac{1}{x+x\ln x}\right)\frac{dx}{\ln x}.
\end{align*}
Since the last integral diverges, the theorem is proved.
\end{pr}

\section{Double points problem for $d=3$}\label{sec4}

We will now focus on the proof of Part (b) of Theorem~\ref{th:OperatorDouble}. Denote
$$
D=\left\{(x,y)\in\RR^3\times\RR^3:|x|>1, |y|>1\right\},
$$
and for $\beta>0$ let
\begin{align*}
J_{\beta}=\iint_{D}&\frac{1}{(1+|x_1+y_1|^{\beta}+|x_2+y_2|^{\beta}+|x_3+y_3|^{\beta})}\\
&\times\frac{dxdy}{(|x_1|^{\alpha_1}+|x_2|^{\alpha_2}+|x_3|^{\alpha_3})(|y_1|^{\alpha_1}+|y_2|^{\alpha_2}+|y_3|^{\alpha_3})}.
\end{align*}
As in Section~\ref{sec3}, we use (\ref{est:Levyexp}) and Theorem~\ref{th:LevyMultiple} to conclude that
\begin{equation}\label{eq:Hausdorff2}
\dimH M_2= 3-\inf\left\{\beta\in(0,3):J_{\beta}<\infty \right\}.
\end{equation}
Hence, Part (b) of Theorem~\ref{th:OperatorDouble} follows from (\ref{eq:Hausdorff2}) and Theorem~\ref{th:Jbeta} below.

\begin{theorem}\label{th:Jbeta}
When $2- \sum_{j=1}^3 \frac 1 {\alpha_j} > 0$,  we have
$$
\inf\left\{\beta>0:J_{\beta}<\infty\right\}=4+\frac{\alpha_1}{\alpha_2}+\frac{\alpha_1}{\alpha_3}-2\alpha_1.
$$
If $2-\sum_{j=1}^3 \frac 1 {\alpha_j} \leq0$, then $J_{\beta}=\infty$ for all $\beta>0$.
\end{theorem}

\begin{pr}
In the proof we assume that $\alpha_3\leq\alpha_2<2$. When $\alpha_3\leq\alpha_2=\alpha_1=2$, the reasoning 
is similar to Part (ii) of the proof of Theorem 4. For $k,n,m\in\NN$, denote
$$
D_{k,n,m}=\left\{(x,y)\in\RR^3\times\RR^3:|x|>1, |y|>1, k-1\leq|x_1+y_1|<k,\right.$$ 
$$\left.n-1\leq|x_2+y_2|<n,m-1\leq|x_3+y_3|<m\right\}.
$$
We have $J_{\beta}<\infty$ if and only if
\begin{align}
\sum_{k=1}^{\infty}\sum_{n=1}^{\infty}\sum_{m=1}^{\infty}\frac{1}{k^{\beta}+n^{\beta}+m^{\beta}}
\iint_{D_{k,n,m}}&\frac{1}{(|x_1|^{\alpha_1}+|x_2|^{\alpha_2}+|x_3|^{\alpha_3})}\label{eq:3series}\\
\times&\frac{dxdy}{(|y_1|^{\alpha_1}+|y_2|^{\alpha_2}+|y_3|^{\alpha_3})}<\infty.\nonumber
\end{align}
Our goal is to prove that if  $2-1/\alpha_1-1/\alpha_2-1/\alpha_3>0$ then
\begin{align}\label{integral}
&\iint_{D_{k,n,m}}\frac{dxdy}{(|x_1|^{\alpha_1}+|x_2|^{\alpha_2}+|x_3|^{\alpha_3})(|y_1|^{\alpha_1}+|y_2|^{\alpha_2}+|y_3|^{\alpha_3})}\\
&\asymp(k^{\alpha_1}+n^{\alpha_2}+m^{\alpha_3})^{1/\alpha_1+1/\alpha_2+1/\alpha_3-2}\nonumber
\end{align}
for sufficiently large $k,n,m$,  whereas the integral in (\ref{integral}) is infinite if $2-1/\alpha_1-1/\alpha_2-1/\alpha_3\leq0$. 
The theorem will then follow from Lemma~\ref{lem:series2} below.

For this purpose, we may and will assume that $|x_1|, |x_2|, |x_3| \geq1$ and $y_1,y_2,y_3\geq1$.
We then have
$$
|y_1-k|\wedge|y_1-k+1|\leq|x_1|\leq y_1+k,
$$
$$
|y_2-n|\wedge|y_2-n+1|\leq|x_2|\leq y_2+n,
$$
$$
|y_3-m|\wedge|y_3-m+1|\leq|x_3|\leq y_3+m.
$$
Furthermore, we may also assume that $|x_1|\geq |y_1-k|$, $|x_2|\geq |y_2-n|$, $|x_3|\geq |y_3-m|$ 
and $|y_1-k|+|y_2-n|+|y_3-m|\geq C>0$.
Therefore, the integral  in (\ref{integral}) can be estimated from below by
\begin{equation}\label{eq:3seriesLow}
\int_1^{\infty}\int_1^{\infty}\int_1^{\infty}\frac{dy_1dy_2dy_3}{(y_1^{\alpha_1}+y_2^{\alpha_2}+y_3^{\alpha_3})((y_1+k)^{\alpha_1}+(y_2+n)^{\alpha_2}+(y_3+m)^{\alpha_3})}
\end{equation}
and from above by
\begin{equation}\label{eq:3seriesUp}
\int_{E_{k,n,m}}\frac{dy}{(y_1^{\alpha_1}+y_2^{\alpha_2}+y_3^{\alpha_3})(|y_1-k|^{\alpha_1}+|y_2-n|^{\alpha_2}+|y_3-m|^{\alpha_3})},
\end{equation}
where
$$
E_{k,n,m}=\left\{y\in\RR^3:\min\left\{y_1,y_2,y_3,|y_1-k|,|y_2-n|,|y_3-m|\right\}\geq1\right\}.
$$
The integral (\ref{eq:3seriesUp}) can be written as
\begin{align*}
&\left(\int_1^{m-1}+\int_{m+1}^{\infty}\right)\left(\int_1^{n-1}+\int_{n+1}^{\infty}\right)\left(\int_1^{k-1}+\int_{k+1}^{\infty}\right)
\frac{1}{y_1^{\alpha_1}+y_2^{\alpha_2}+y_3^{\alpha_3}}\\
&\qquad \times\frac{dy_1dy_2dy_3}{|y_1-k|^{\alpha_1}+|y_2-n|^{\alpha_2}+|y_3-m|^{\alpha_3}},
\end{align*}
which is equal to the sum of the following integrals
\begin{itemize}
\item[1)]$\displaystyle\int_1^{m-1}\int_1^{n-1}\int_1^{k-1}$,
\item[2)]$\displaystyle\int_1^{m-1}\int_{n+1}^{\infty}\int_{k+1}^{\infty},\quad\int_{m+1}^{\infty}\int_1^{n-1}\int_{k+1}^{\infty},
\quad\int_{m+1}^{\infty}\int_{n+1}^{\infty}\int_1^{k-1}$,
\item[3)]$\displaystyle\int_1^{m-1}\int_1^{n-1}\int_{k+1}^{\infty},\quad\int_1^{m-1}\int_{n+1}^{\infty}\int_1^{k-1},
\quad\int_{m+1}^{\infty}\int_1^{n-1}\int_1^{k-1}$,
\item[4)] $\displaystyle\int_{m+1}^{\infty}\int_{n+1}^{\infty}\int_{k+1}^{\infty}$.
\end{itemize}
Since the integral in 4) is the same as the one in the lower bound (\ref{eq:3seriesLow}), 
we start by establishing desired upper and lower bounds  as in (\ref{integral}) for the integral 
in 4). To simplify the notation, denote $\eta:=k^{\alpha_1}+n^{\alpha_2}+m^{\alpha_3}$. 
One can verify that
\begin{align*}
&\int_{m+1}^{\infty}\int_{n+1}^{\infty}\int_{k+1}^{\infty}
\frac{dy_1dy_2dy_3}{(y_1^{\alpha_1}+y_2^{\alpha_2}+y_3^{\alpha_3})(|y_1-k|^{\alpha_1}+|y_2-n|^{\alpha_2}+|y_3-m|^{\alpha_3})}\\
&\asymp\int_{1}^{\infty}\int_{1}^{\infty}\int_{1}^{\infty}
\frac{dx_1dx_2dx_3}{(x_1^{\alpha_1}+x_2^{\alpha_2}+x_3^{\alpha_3})(x_1^{\alpha_1}+x_2^{\alpha_2}+x_3^{\alpha_3}+\eta)}\\
&\asymp\int_{1}^{\infty}\int_{1}^{\infty}\int_{1}^{\infty}\left(\frac{1}{x_1^{\alpha_1}}\wedge\frac{1}{x_2^{\alpha_2}} \wedge\frac{1}{x_3^{\alpha_3}}\right)
\left(\frac{1}{x_1^{\alpha_1}}\wedge\frac{1}{x_2^{\alpha_2}} \wedge\frac{1}{x_3^{\alpha_3}} \wedge\frac{1}{\eta}\right)\, dx_1dx_2dx_3\\
&=\int_{1}^{\infty}\int_{1}^{x_3^{\alpha_3/\alpha_2}}\int_{1}^{x_2^{\alpha_2/\alpha_1}}\frac{1}{x_3^{\alpha_3}}
\left[\frac{1}{x_3^{\alpha_3}}\wedge\frac{1}{\eta}\right]dx_1dx_2dx_3\\
&\qquad +\int_{1}^{\infty}\int_{x_3^{\alpha_3/\alpha_2}}^{\infty}\int_{1}^{x_2^{\alpha_2/\alpha_1}}\frac{1}{x_2^{\alpha_2}}
\left[\frac{1}{x_2^{\alpha_2}}\wedge\frac{1}{\eta}\right]dx_1dx_2dx_3\\
&\qquad +\int_{1}^{\infty}\int_{x_2^{\alpha_2/\alpha_1}}^{\infty}\int_{1}^{x_1^{\alpha_1/\alpha_3}}\frac{1}{x_1^{\alpha_1}}
\left[\frac{1}{x_1^{\alpha_1}}\wedge\frac{1}{\eta}\right]dx_3dx_1dx_2\\
&\qquad +\int_{1}^{\infty}\int_{x_2^{\alpha_2/\alpha_1}}^{\infty}\int_{x_1^{\alpha_1/\alpha_3}}^{\infty}\frac{1}{x_3^{\alpha_3}}
\left[\frac{1}{x_3^{\alpha_3}}\wedge\frac{1}{\eta}\right]dx_3dx_1dx_2=:I_1+I_2+I_3+I_4.
\end{align*}
For $I_1$, by breaking the integral according to $x_3 \le \eta^{1/\alpha_3}$ and $x_3 > \eta^{1/\alpha_3}$, we can verify that 
$I_1$ is convergent if and only if $2-1/\alpha_1-1/\alpha_2-1/\alpha_3>0$ and in the later case, 
$$
I_1 \asymp\eta^{1/\alpha_1+1/\alpha_2+1/\alpha_3-2}.
$$
This also proves the second part of the theorem. 

Next we assume $2-1/\alpha_1-1/\alpha_2-1/\alpha_3>0$.
Then it is elementary to verify that
$$
I_2,\, I_3 ,\, I_4 \leq C\eta^{1/\alpha_1+1/\alpha_2+1/\alpha_3-2}.
$$
Hence the integral  in 4) satisfies the same bounds.

Next, consider the integral in 1). Noticing that
$$
\int_1^{m-1}\int_1^{n-1}\int_1^{k/2} (\cdots) =\int_1^{m-1}\int_1^{n-1}\int_{k/2}^{k-1} (\cdots),
$$
we have  
\begin{align*}
&\int_1^{m-1}\int_1^{n-1}\int_1^{k-1}\frac{dy_1dy_2dy_3}{(y_1^{\alpha_1}+y_2^{\alpha_2}+y_3^{\alpha_3})
(|y_1-k|^{\alpha_1}+|y_2-n|^{\alpha_2}+|y_3-m|^{\alpha_3})}\\
&\leq C\int_1^{k/2}\frac{dy_1}{(k-y_1)^{\alpha_1}}\int_1^{n-1}\int_1^{m-1}
\frac{dy_3dy_2 }{(y_1^{\alpha_1/\alpha_3}+y_2^{\alpha_2/\alpha_3}+y_3)^{\alpha_3}}.
\end{align*}
Integrating out $dy_3$ and $dy_2$ we see that the last integral is at most 
$$
C k^{-\alpha_1} \int_1^{k/2} {(n+y_1^{\alpha_1/\alpha_2})^{1+\alpha_2/\alpha_3-\alpha_2}}\,dy_1
\leq C k^{-\alpha_1}(k^{\alpha_1}+n^{\alpha_2})^{1/\alpha_1+1/\alpha_2+1/\alpha_3-1}.
$$
By symmetry, the integral in 1) is also less than
$$
Cn^{-\alpha_2}(k^{\alpha_1}+n^{\alpha_2})^{1/\alpha_1+1/\alpha_2+1/\alpha_3-1}.
$$
Combining two terms we see that the integral in 1) is at most
\begin{align*}
& C(k^{-\alpha_1}\wedge n^{-\alpha_2})(k^{\alpha_1}+n^{\alpha_2})^{1/\alpha_1+1/\alpha_2+1/\alpha_3-1}\\
&\leq C(k^{\alpha_1}+n^{\alpha_2})^{1/\alpha_1+1/\alpha_2+1/\alpha_3-2}.
\end{align*}
Since similar estimates work with pairs $k^{\alpha_1},m^{\alpha_3}$ and $n^{\alpha_2},m^{\alpha_3}$, we obtain
the following majorant for the integral in 1):
\begin{align*}
&C\min\left\{k^{\alpha_1}+n^{\alpha_2},\, (k^{\alpha_1}+ m^{\alpha_3}),\, (n^{\alpha_2}+m^{\alpha_3})\right\}^{1/\alpha_1+1/\alpha_2+1/\alpha_3-2}\\
&\leq C(k^{\alpha_1}+n^{\alpha_2}+m^{\alpha_3})^{1/\alpha_1+1/\alpha_2+1/\alpha_3-2}.
\end{align*}
Next, by the symmetry, it is enough to consider only one integral of type 2) and one of type 3).
Consider first an integral of type 3) as follows.
\begin{equation}\label{inttype3}
\int_1^{m-1}\int_1^{n-1}\int_{k+1}^{\infty}\frac{dy_1dy_2dy_3}{(y_1^{\alpha_1}+y_2^{\alpha_2}+y_3^{\alpha_3})(|y_1-k|^{\alpha_1}+|y_2-n|^{\alpha_2}+|y_3-m|^{\alpha_3})}
\end{equation}
$$
=\left(\int_1^{m/2}\int_1^{n/2}\int_{1}^{\infty}+\int_{m/2}^{m-1}\int_1^{n/2}\int_{1}^{\infty}+\int_1^{m/2}\int_{n/2}^{n-1}\int_{1}^{\infty}+\int_{m/2}^{m-1}\int_{n/2}^{n-1}\int_{1}^{\infty}\right)
$$
$$
\times\frac{dy_1dy_2dy_3}{(y_1^{\alpha_1}+y_2^{\alpha_2}+y_3^{\alpha_3})((y_1+k)^{\alpha_1}+(n-y_2)^{\alpha_2}+(m-y_3)^{\alpha_3})}
=:J_1+J_2+J_3+J_4.
$$
We obtain
\begin{align*}
J_1&=\int_1^{m/2}\int_1^{n/2}\int_{1}^{\infty}\frac{dy_1dy_2dy_3}{(y_1^{\alpha_1}+y_2^{\alpha_2}+y_3^{\alpha_3})
((y_1+k)^{\alpha_1}+(n-y_2)^{\alpha_2}+(m-y_3)^{\alpha_3})}\\
&\leq C\int_1^{m/2}\int_1^{n/2}\int_{1}^{\infty}\frac{dy_1dy_2dy_3}{(y_2^{\alpha_2}+y_3^{\alpha_3})(y_1^{\alpha_1}+m^{\alpha_3})}\\
&\leq C\int_{1}^{\infty}\frac{dy_1}{(y_1+m^{\alpha_3/\alpha_1})^{\alpha_1}}\int_1^{m/2}\int_1^{n/2}\frac{dy_2dy_3}{(y_2+y_3^{\alpha_3/\alpha_2})^{\alpha_2}}\\
&\leq C \left(m^{\alpha_3}\right)^{1/\alpha_1+1/\alpha_2+1/\alpha_3-2}.
\end{align*}
For $J_2$, a simple change of variable yields
\begin{align*}
J_2=&\int_1^{m/2}\int_1^{n/2}\int_{1}^{\infty}\frac{dy_1dy_2dr}{(y_1^{\alpha_1}+y_2^{\alpha_2}+(m-r)^{\alpha_3})((y_1+k)^{\alpha_1}+(n-y_2)^{\alpha_2}+r^{\alpha_3})}\\
\leq& C\int_1^{m/2}\int_1^{n/2}\int_{1}^{\infty}\frac{dy_1dy_2dr}{(y_2^{\alpha_2}+m^{\alpha_3})(y_1^{\alpha_1}+r^{\alpha_3})}\\
\leq& C \left(m^{\alpha_3}\right)^{1/\alpha_1+1/\alpha_2+1/\alpha_3-2}.
\end{align*}
In a similar manner we obtain
\begin{align*}
J_3&=\int_1^{m/2}\int_1^{n/2}\int_{1}^{\infty}\frac{dy_1dsdy_3}{(y_1^{\alpha_1}+(n-s)^{\alpha_2}+y_3^{\alpha_3})((y_1+k)^{\alpha_1}+s^{\alpha_2}+(m-y_3)^{\alpha_3})}\\
&\leq C\int_1^{m/2}\int_1^{n/2}\int_{1}^{\infty}\frac{dy_1dsdy_3}{(y_1^{\alpha_1}+y_3^{\alpha_3})(s^{\alpha_2}+m^{\alpha_3})},
\end{align*}
so the estimate is the same as for $J_2$. Finally,
\begin{align*}
J_4 &=\int_1^{m/2}\int_1^{n/2}\int_{1}^{\infty}\frac{dy_1dsdr}{(y_1^{\alpha_1}+(n-s)^{\alpha_2}+(m-r)^{\alpha_3})((y_1+k)^{\alpha_1}+s^{\alpha_2}+r^{\alpha_3})}\\
&\leq C\int_1^{m/2}\int_1^{n/2}\int_{1}^{\infty}\frac{dy_1dsdr}{(y_1^{\alpha_1}+m^{\alpha_3})(s^{\alpha_2}+r^{\alpha_3})},
\end{align*}
which is the same integral appeared in the estimation of $J_1$. Hence we have proved that the integral (\ref{inttype3}) is less than 
$C \left(m^{\alpha_3}\right)^{1/\alpha_1+1/\alpha_2+1/\alpha_3-2}$. By symmetry, it is also less than
$C \left(n^{\alpha_2}\right)^{1/\alpha_1+1/\alpha_2+1/\alpha_3-2}$. 

In order to obtain a similar upper bound 
in terms of  $k^{\alpha_1}$ instead of $m^{\alpha_3}$ we observe that
\begin{align*}
&\int_1^{m-1}\int_1^{n-1}\int_{k+1}^{\infty}\frac{dy_1dy_2dy_3}{(y_1^{\alpha_1}+y_2^{\alpha_2}+y_3^{\alpha_3})
(|y_1-k|^{\alpha_1}+|y_2-n|^{\alpha_2}+|y_3-m|^{\alpha_3})}\\
&\leq n\int_1^{m-1}\int_1^{\infty}\frac{dy_1dy_3}{(y_1^{\alpha_1}+y_3^{\alpha_3})((y_1+k)^{\alpha_1}+(m-y_3)^{\alpha_3})}.
\end{align*}
The last double integral is of the same type as the integral (\ref{eq:int3}) in the proof of Theorem~\ref{th:Ibeta}, and 
therefore, it is less than
$$
Cn\left(k^{\alpha_1}\right)^{1/\alpha_1+1/\alpha_3-2}=C\left(\frac{n^{\alpha_2}}{k^{\alpha_1}}\right)^{1/\alpha_2}\left(k^{\alpha_1}\right)^{1/\alpha_1+1/\alpha_2+1/\alpha_3-2}.
$$
If $k^{\alpha_1}>n^{\alpha_2}$, then we get the desired estimate. On the other hand, when $k^{\alpha_1}\leq n^{\alpha_2}$, then
$$
\left(n^{\alpha_2}\right)^{1/\alpha_1+1/\alpha_2+1/\alpha_3-2}\leq\left(k^{\alpha_1}\right)^{1/\alpha_1+1/\alpha_2+1/\alpha_3-2},
$$
since $1/\alpha_1+1/\alpha_2+1/\alpha_3-2<0$ by our assumption, and the upper bound follows from 
the previous part of the proof. Therefore, the minimum of obtained upper bounds gives
$C(k^{\alpha_1}+n^{\alpha_2}+m^{\alpha_3})^{1/\alpha_1+1/\alpha_2+1/\alpha_3-2}$ as 
a majorant for the integrals of type 3).

Finally, we consider an integral of type 2). After a change of variables we have
$$
\int_{m+1}^{\infty}\int_{n+1}^{\infty}\int_1^{k-1}\frac{dy_1dy_2dy_3}{(y_1^{\alpha_1}+y_2^{\alpha_2}+y_3^{\alpha_3})
(|y_1-k|^{\alpha_1}+|y_2-n|^{\alpha_2}+|y_3-m|^{\alpha_3})}
$$
$$
=\int_1^{\infty}\int_1^{\infty}\int_1^{k-1}\frac{dtdsdr}{(t^{\alpha_1}+s^{\alpha_2}+r^{\alpha_3})((k-t)^{\alpha_1}+(n+s)^{\alpha_2}+(m+r)^{\alpha_3})}.
$$
Furthermore, by breaking the integration interval $[1, k-1]$ in $dt$ into $[1, k/2]$ and $[k/2, k-1]$, we derive
\begin{align*}
&\int_1^{\infty}\int_1^{\infty}\int_1^{k-1}\frac{dtdsdr}{(t^{\alpha_1}+s^{\alpha_2}+r^{\alpha_3})((k-t)^{\alpha_1}+(n+s)^{\alpha_2}+(m+r)^{\alpha_3})}\\
&\leq C\int_1^{\infty}\int_1^{\infty}\int_1^{k/2}\frac{dtdsdr}{(t^{\alpha_1}+s^{\alpha_2}+r^{\alpha_3})(k^{\alpha_1}+s^{\alpha_2}+r^{\alpha_3})}.
\end{align*}
One can show that the last term is less than $C(k^{\alpha_1})^{1/\alpha_1+1/\alpha_2+1/\alpha_3-2}$, 
which implies the same estimate for the initial integral of type 2). Since the method is similar to the case of 
the integral of type 4), we omit the details. Therefore, given Lemma \ref{lem:series2}, the proof of Theorem \ref{th:Jbeta} 
is finished.
\end{pr}

In order to prove Lemma \ref{lem:series2} we will make use of the following inequality, whose proof is elementary and is omitted.
\begin{lemma}\label{lem:alpha}
For $2\geq\alpha_1\geq\alpha_2\geq\alpha_3>1$ we have
$$
\max\left\{2+\frac{2\alpha_3}{\alpha_1}+\frac{2\alpha_3}{\alpha_2}-2\alpha_3,3+\frac{2\alpha_2}
{\alpha_1}+\frac{\alpha_2}{\alpha_3}-2\alpha_2\right\}
\leq 4+\frac{\alpha_1}{\alpha_2}+\frac{\alpha_1}{\alpha_3}-2\alpha_1.
$$
\end{lemma}

\begin{lemma}\label{lem:series2}
Let $\gamma:=2-1/\alpha_1-1/\alpha_2-1/\alpha_3$.
We have
$$
\inf\left\{\beta\in(0,3):\sum_{k=1}^{\infty}\sum_{n=1}^{\infty}\sum_{m=1}^{\infty}\frac{(k^{\alpha_1}
+ n^{\alpha_2}+m^{\alpha_3})^{-\gamma}}{k^{\beta}+n^{\beta}+m^{\beta}}<\infty\right\}
=4+\frac{\alpha_1}{\alpha_2}+\frac{\alpha_1}{\alpha_3}-2\alpha_1.
$$
\end{lemma}

\begin{pr}
The convergence of the series is equivalent to the convergence of the integral
\begin{align*}
&\int_1^{\infty}\int_1^{\infty}\int_1^{\infty}\left[ x\vee y \vee z\right]^{-\beta}\left[ x^{\alpha_1}\vee y^{\alpha_2} \vee z^{\alpha_3}\right]^{-\gamma}dxdydz\\
&=\int_1^{\infty}\int_1^{\infty}\int_1^{y^{\alpha_2/\alpha_1}}(y\vee z)^{-\beta}(y^{\alpha_2}\vee z^{\alpha_3})^{-\gamma}dxdydz\\
&+\int_1^{\infty}\int_1^{\infty}\int_{y^{\alpha_2/\alpha_1}}^y(y\vee z)^{-\beta}(x^{\alpha_1}\vee z^{\alpha_3})^{-\gamma}dxdydz\\
&+\int_1^{\infty}\int_1^{\infty}\int_y^{\infty}(x\vee z)^{-\beta}(x^{\alpha_1}\vee z^{\alpha_3})^{-\gamma}dxdydz=:I_1+I_2+I_3.
\end{align*}
We have $I_1<\infty$ if and only if
\begin{align*}
\infty&>\int_1^{\infty}\int_1^{\infty}y^{\alpha_2/\alpha_1}(y\vee z)^{-\beta}(y^{\alpha_2}\vee z^{\alpha_3})^{-\gamma}dydz\\
&=\int_1^{\infty}\int_1^{z^{\alpha_3/\alpha_2}}y^{\alpha_2/\alpha_1}z^{-\beta-\alpha_3\gamma}dydz+
\int_1^{\infty}\int_{z^{\alpha_3/\alpha_2}}^zy^{\alpha_2/\alpha_1-\alpha_2\gamma}z^{-\beta}dydz\\
&+\int_1^{\infty}\int_z^{\infty}y^{\alpha_2/\alpha_1-\alpha_2\gamma-\beta}dydz=:I_1^{(1)}+I_1^{(2)}+I_1^{(3)}.
\end{align*}
First, it can be verified that  $I_1^{(1)}<\infty$ if and only if
$\alpha_3/\alpha_2(\alpha_2/\alpha_1+1)-\beta-\alpha_3\gamma<-1$. This gives the inequality $\beta>2+2\alpha_3/\alpha_1+2\alpha_3/\alpha_2-2\alpha_3$.
Furthermore, a necessary condition for $I_1^{(3)}<\infty$ is $\alpha_2/\alpha_1-\alpha_2\gamma-\beta<-1$. Assuming this we get
$$
I_1^{(3)}=\int_1^{\infty}\frac{z^{\alpha_2/\alpha_1-\alpha_2\gamma-\beta+1}}{\alpha_2\gamma+\beta-\alpha_2/\alpha_1-1}dz.
$$
Hence, $I_1^{(3)}<\infty$ if and only if $\alpha_2/\alpha_1-\alpha_2\gamma-\beta+1<-1$, which gives $\beta>3+2\alpha_2/\alpha_1+\alpha_2/\alpha_3-2\alpha_2$.

In order to estimate $I_1^{(2)}$, we observe that $\alpha_2/\alpha_1-\alpha_2\gamma>-1$. Hence
$$
I_1^{(2)} \asymp \int_1^{\infty}\frac{z^{1+\alpha_2/\alpha_1-\alpha_2\gamma-\beta}}{\alpha_2/\alpha_1-\alpha_2\gamma+1}dz,
$$
and the last integral is finite if and only if $\beta>3+2\alpha_2/\alpha_1+\alpha_2/\alpha_3-2\alpha_2$. Therefore, we have proved that $I_1<\infty$ if and only if
$$
\beta>\max\left\{2+\frac{2\alpha_3}{\alpha_1}+\frac{2\alpha_3}{\alpha_2}-2\alpha_3,3+\frac{2\alpha_2}{\alpha_1}+\frac{\alpha_2}{\alpha_3}-2\alpha_2\right\}.
$$
Next we rewrite $I_2$ as
\begin{align*}
I_2&=\int_1^{\infty}\int_1^{z^{\alpha_3/\alpha_1}}\int_{y^{\alpha_2/\alpha_1}}^yz^{-\beta-\gamma\alpha_3}dxdydz+
\int_1^{\infty}\int_{z^{\alpha_3/\alpha_1}}^{z^{\alpha_3/\alpha_2}}\int_{y^{\alpha_2/\alpha_1}}^{z^{\alpha_3/\alpha_1}}z^{-\beta-\gamma\alpha_3}dxdydz\\
&+\int_1^{\infty}\int_{z^{\alpha_3/\alpha_1}}^{z^{\alpha_3/\alpha_2}}\int_{z^{\alpha_3/\alpha_1}}^yz^{-\beta}x^{-\gamma\alpha_1}dxdydz
+\int_1^{\infty}\int_{z^{\alpha_3/\alpha_2}}^{z}\int_{y^{\alpha_2/\alpha_1}}^yz^{-\beta}x^{-\gamma\alpha_1}dxdydz\\
&+\int_1^{\infty}\int_{z}^{\infty}\int_{y^{\alpha_2/\alpha_1}}^yy^{-\beta}x^{-\gamma\alpha_1}dxdydz=:I_2^{(1)}+I_2^{(2)}+I_2^{(3)}+I_2^{(4)}+I_2^{(5)}.
\end{align*}
We get
$$
I_2^{(1)}\leq\int_1^{\infty}\int_1^{z^{\alpha_3/\alpha_1}}yz^{-\beta-\gamma\alpha_3}dydz
\leq\int_1^{\infty}z^{-\beta-\gamma\alpha_3+2\alpha_3/\alpha_1}dz,
$$
and the last integral is finite if and only if $\beta>2+3\alpha_3/\alpha_1+\alpha_3/\alpha_2+2\alpha_3$. 
However, since $\alpha_3/\alpha_1\leq\alpha_3/\alpha_2$, the condition
$\beta>2+2\alpha_3/\alpha_1+2\alpha_3/\alpha_2+2\alpha_3$ implies $I_2^{(1)}<\infty$. Secondly,
$$
I_2^{(2)}\leq\int_1^{\infty}\int_{z^{\alpha_3/\alpha_1}}^{z^{\alpha_3/\alpha_2}}z^{-\beta-\gamma\alpha_3+\alpha_3/\alpha_1}dydz
\leq\int_1^{\infty}z^{-\beta-\gamma\alpha_3+\alpha_3/\alpha_1+\alpha_3/\alpha_2}dz.
$$
Hence the condition $\beta>2+2\alpha_3/\alpha_1+2\alpha_3/\alpha_2+2\alpha_3$ implies $I_2^{(2)}<\infty$. Next we have
$$
I_2^{(3)}\vee I_2^{(4)}\leq \int_1^{\infty}\int_1^z\int_1^y z^{-\beta}x^{-\gamma\alpha_1}dxdydz.
$$
When $\alpha_3<2$, then $-\gamma\alpha_1>-1$ and the last term is less than
$$
 C\int_1^{\infty}\int_1^z z^{-\beta}y^{1-\gamma\alpha_1}dydz
 \leq C\int_1^{\infty} z^{2-\beta-\gamma\alpha_1}dz.
$$
The last integral is finite if and only if $2-\beta-\gamma\alpha_1<-1$, which is equivalent to $\beta>4+{\alpha_1}/{\alpha_2}+{\alpha_1}/{\alpha_3}-2\alpha_1$.
If $\alpha_3=2$, then
$$
\int_1^{\infty}\int_1^z\int_1^y z^{-\beta}x^{-\gamma\alpha_1}dxdydz=\int_1^{\infty}\int_1^z z^{-\beta}\ln(y)dydz\leq\int_1^{\infty} z^{1-\beta}\ln(z)dz,
$$
which again gives the condition $\beta>2=4+{\alpha_1}/{\alpha_2}+{\alpha_1}/{\alpha_3}-2\alpha_1$.
This then implies $I_2^{(3)}, I_2^{(4)}<\infty$. Finally,
$$
I_2^{(5)}\leq C\int_1^{\infty}\int_z^{\infty} y^{1-\beta-\gamma\alpha_1}dydz.
$$
For $1-\beta-\gamma\alpha_1<-1$ the last therm is less than
$$
C\int_1^{\infty}z^{2-\beta-\gamma\alpha_1}dz.
$$
Therefore, the inequality $\beta>4+{\alpha_1}/{\alpha_2}+{\alpha_1}/{\alpha_3}-2\alpha_1$ 
implies $I_2^{(5)}<\infty$. Hence, we have proved  $I_2<\infty$ provided
$$
\beta>\max\left\{2+\frac{2\alpha_3}{\alpha_1}+\frac{2\alpha_3}{\alpha_2}-2\alpha_3,
4+\frac{\alpha_1}{\alpha_2}+\frac{\alpha_1}{\alpha_3}-2\alpha_1\right\}.
$$

Consider $I_3$. We have
\begin{align*}
I_3&=\int_1^{\infty}\int_1^{z^{\alpha_3/\alpha_1}}\int_y^{z^{\alpha_3/\alpha_1}}z^{-\beta-\gamma\alpha_3}dxdydz
+\int_1^{\infty}\int_1^{z^{\alpha_3/\alpha_1}}\int_{z^{\alpha_3/\alpha_1}}^zz^{-\beta}x^{-\gamma\alpha_1}dxdydz\\
&\qquad+\int_1^{\infty}\int_1^{z^{\alpha_3/\alpha_1}}\int_z^{\infty}x^{-\beta-\gamma\alpha_1}dxdydz
+\int_1^{\infty}\int_{z^{\alpha_3/\alpha_1}}^z\int_y^zz^{-\beta}x^{-\gamma\alpha_1}dxdydz\\
&\qquad +\int_1^{\infty}\int_{z^{\alpha_3/\alpha_1}}^z\int_z^{\infty}x^{-\beta-\gamma\alpha_1}dxdydz
+\int_1^{\infty}\int_z^{\infty}\int_y^{\infty}x^{-\beta-\gamma\alpha_1}dxdydz\\
&=:I_3^{(1)}+I_3^{(2)}+I_3^{(3)}+I_3^{(4)}+I_3^{(5)}+I_3^{(6)}.
\end{align*}
We obtain
$$
I_3^{(1)}\leq\int_1^{\infty}\int_1^{z^{\alpha_3/\alpha_1}}z^{-\beta-\gamma\alpha_3+\alpha_3/\alpha_1}dydz
\leq\int_1^{\infty}z^{-\beta-\gamma\alpha_3+2\alpha_3/\alpha_1}dz.
$$
The last integral is finite if and only if $-\beta-\gamma\alpha_3+2\alpha_3/\alpha_1<-1$ which gives
$\beta>2+3\alpha_3/\alpha_1+\alpha_3/\alpha_2-2\alpha_3$. Since
$$
2+3\alpha_3/\alpha_1+\alpha_3/\alpha_2-2\alpha_3<2+2\alpha_3/\alpha_1+2\alpha_3/\alpha_2-2\alpha_3,
$$
the condition $\beta>2+2\alpha_3/\alpha_1+2\alpha_3/\alpha_2-2\alpha_3$ implies $I_3^{(1)}<\infty$. 
Next, when $\alpha_3<2$ then we get $-\gamma\alpha_1>-1$ and
$$
I_3^{(2)}\leq C\int_1^{\infty}\int_1^{z^{\alpha_3/\alpha_1}}z^{1-\beta-\gamma\alpha_1}dydz
\leq C\int_1^{\infty}z^{1+\alpha_3/\alpha_1-\beta-\gamma\alpha_1}dz.
$$
The last integral is finite if and only if $1+\alpha_3/\alpha_1-\beta-\gamma\alpha_1<-1$, which gives 
$\beta>3+\alpha_3/\alpha_1+\alpha_1/\alpha_2+\alpha_1/\alpha_3-2\alpha_1$. 

Proceeding in a similar manner we can show that $\beta>4+\alpha_1/\alpha_2+\alpha_1/\alpha_3-2\alpha_1$ 
implies  $I_3^{(3)},\,  I_3^{(4)},\, I_3^{(5)}<\infty$.

Finally,
$$
I_3^{(6)}=\int_1^{\infty}\int_z^{\infty}\frac{y^{1-\beta-\gamma\alpha_1}}{\beta+\gamma\alpha_1-1}dydz
=\int_1^{\infty}\frac{z^{2-\beta-\gamma\alpha_1}}{(\beta+\gamma\alpha_1-1)(\beta+\gamma\alpha_1-2)}dz.
$$
Hence, $I_3^{(6)}<\infty$ if and only if $2-\beta-\gamma\alpha_1<-1$, which gives 
$\beta>4+\alpha_1/\alpha_2+\alpha_1/\alpha_3-2\alpha_1$.

In summary, we have proved the following
$$
I_1<\infty\Leftrightarrow\beta>\max\left\{2+\frac{2\alpha_3}{\alpha_1}+\frac{2\alpha_3}
{\alpha_2}-2\alpha_3,\, 3+\frac{2\alpha_2}{\alpha_1}+\frac{\alpha_2}{\alpha_3}-2\alpha_2\right\},
$$
$$
\beta>\max\left\{2+\frac{2\alpha_3}{\alpha_1}+\frac{2\alpha_3}{\alpha_2}-2\alpha_3,\,
4+\frac{\alpha_1}{\alpha_2}+\frac{\alpha_1}{\alpha_3}-2\alpha_1\right\}\Rightarrow I_2,I_3<\infty,
$$
and $I_3<\infty\Rightarrow I_3^{(6)}<\infty\Rightarrow\beta>4+\alpha_1/\alpha_2+\alpha_1/\alpha_3-2\alpha_1$. 
The final conclusion follows from Lemma~\ref{lem:alpha}.
\end{pr}

Our last result gives a necessary and sufficient condition for the existence of double points of $X$ in $\RR^3$. 
It differs slightly from the case $d=2$, i.e., $X$ may possess double points even if the Hausdorff dimension of 
$M_2$ is 0.

According to the decomposition described in Section~\ref{sec2}, the stability exponent of $X$ satisfies $B=PDP^{-1}$. 
In the present case, the matrix $D$ can have the following forms
\begin{enumerate}
\item[(a)] $\left(
	\begin{matrix}
		a_1 & 0 & 0 \\
		0 & a_2 & 0 \\
		0 & 0   & a_3 \\
	\end{matrix}
	\right); $ \\\\
\item[(b)] $\left(
	\begin{matrix}
		a_1 & 0 & 0 \\
		1 & a_1 & 0 \\
		0 & 0   & a_2 \\
	\end{matrix}
	\right) $ or $\left(
	\begin{matrix}
		a_1 & -b_1 & 0 \\
		b_1 & a_1 & 0 \\
		0 & 0   & a_2 \\
	\end{matrix}
	\right); $\\\\
\item[(c)] $\left(
	\begin{matrix}
		a_1 & 0 & 0 \\
		0 & a_2 & 0 \\
		0 & 1   & a_2 \\
	\end{matrix}
	\right) $ or $\left(
	\begin{matrix}
		a_1 & 0 & 0 \\
		0 & a_2 & -b_2 \\
		0 & b_2   & a_2 \\
	\end{matrix}
	\right); $\\\\
	\item[(d)] $\left(
	\begin{matrix}
		a & 0 & 0 \\
		1 & a & 0 \\
		0 & 1   & a \\
	\end{matrix}
	\right). $
\end{enumerate}

\begin{theorem}\label{th:existd3}
The existence of double points of $X$ depends on the cases (a)-(d) as follows:
\begin{itemize}
\item In Cases (a), (b) and (c), $M_2 \ne \emptyset$  if and only if $2-1/\alpha_1-1/\alpha_2-1/\alpha_3>0$. 
\item In  Case (d), $M_2 \ne \emptyset$  if and only if $\alpha_1\geq 3/2$.
\end{itemize}
\end{theorem}

\begin{pr}
By Theorem \ref{th:MultipleExistence} we have $M_2\neq\emptyset$ if and only if
\begin{equation}\label{eq:intPsi2}
\int_{\RR^3}\left(\frac{1}{1+\Psi(\xi)}\right)^2d\xi<\infty.
\end{equation}
By \cite[(4.9), (4.14), (4.15), (4.16)]{MX}, we have the following estimates of the L\'evy 
exponent $\Psi(\xi)$ when $\|\xi\|\to\infty$, 
depending on Cases (a)-(d):
\begin{enumerate}
\item[(a)] $\alpha_1=1/a_1$, $\alpha_2=1/a_2$, $\alpha_3=1/a_3$, and
$$\displaystyle\Psi(\xi)\asymp |\xi_1|^{\alpha_1}+|\xi_2|^{\alpha_2}+|\xi_3|^{\alpha_3};$$
\item[(b)] $\alpha_1=\alpha_2=1/a_1$, $\alpha_3=1/a_2$, and
$$\displaystyle\Psi(\xi)\asymp |\xi_1|^{\alpha_1}+|\xi_2|^{\alpha_1}(\ln\|\xi\|)^{\alpha_1}+|\xi_3|^{\alpha_2};$$
\item[(c)] $\alpha_1=1/a_1$, $\alpha_2=\alpha_3=1/a_2$, and
$$\displaystyle\Psi(\xi)\asymp |\xi_1|^{\alpha_1}+|\xi_2|^{\alpha_2}+|\xi_3|^{\alpha_2}(\ln\|\xi\|)^{\alpha_2};$$
\item[(d)] $\alpha_1=\alpha_2=\alpha_3=1/a$, and
$$\displaystyle\Psi(\xi)\asymp |\xi_1|^{\alpha_1}+|\xi_2|^{\alpha_1}(\ln\|\xi\|)^{\alpha_1}+|\xi_3|^{\alpha_1}(\ln\|\xi\|)^{2\alpha_1}.$$
\end{enumerate}
In Case (a), it follows from the proof of Theorem~\ref{th:Jbeta} and the estimates of the integral (\ref{eq:3seriesLow}) 
that (\ref{eq:intPsi2}) holds if and only if
$2-1/\alpha_1-1/\alpha_2-1/\alpha_3>0$. This proves the theorem for Case (a).

For Case (d),  (\ref{eq:intPsi2}) holds if and only if
$$
\int_1^{\infty}\int_1^{\infty}\int_1^{\infty}\frac{dxdydz}{\left(x+y\ln\sqrt{x^2+y^2+z^2}+z\ln^2\sqrt{x^2+y^2+z^2}\right)^{2\alpha_1}}<\infty,
$$
which, in turn, is equivalent to  
\begin{equation}\label{eq:intlog1}
\begin{split}
&\int_1^{\infty}\int_1^{x}\int_1^{x}\frac{dzdydx}{\left(x+y\ln x+z\ln^2 x\right)^{2\alpha_1}}
+\int_1^{\infty}\int_1^{y}\int_1^{y}\frac{dzdxdy}{\left(x+y\ln y+z\ln^2 y\right)^{2\alpha_1}}\\
&\qquad \qquad 
+\int_1^{\infty}\int_1^{z}\int_1^{z}\frac{dxdydz}{\left(x+y\ln z+z\ln^2 z\right)^{2\alpha_1}}< \infty.
\end{split}
\end{equation}
For the first integral in (\ref{eq:intlog1}), we have 
\begin{align*}
\int_1^{\infty}\int_1^{x}\int_1^{x}\frac{dzdydx}{\left(x+y\ln x+z\ln^2 x\right)^{2\alpha_1}}
&\geq\int_1^{\infty}\int_1^{x}\int_1^{y}\frac{dzdydx}{\left(x+y\ln x+z\ln^2 x\right)^{2\alpha_1}}\\
&=\sum_{n=0}^{\infty}\int_{2^n}^{2^{n+1}}dx\int_1^{x}\int_1^{y}\frac{dzdy}{\left(x+y\ln x+z\ln^2 x\right)^{2\alpha_1}}\\
&\geq \sum_{n=0}^{\infty}\int_1^{2^n}\int_1^{y}\frac{C2^ndzdy}{\left(2^{n}+ny+n^2z\right)^{2\alpha_1}}.
\end{align*}
When $\alpha_1<3/2$, one can verify that the last series, thus (\ref{eq:intlog1}), is infinite. This proves that 
$\alpha_1\geq3/2$ is a necessary condition for (\ref{eq:intPsi2}) in the case (d). 

To prove sufficiency, it is enough to show that the three integrals in (\ref{eq:intlog1}) are finite
for $\alpha_1=3/2$. Since the method is similar, we only consider the first integral.  
\begin{align*}
\int_1^{\infty}\int_1^{x}\int_1^{x}\frac{dzdydx}{\left(x+y\ln x+z\ln^2 x\right)^3}
&=\sum_{n=0}^{\infty}\int_{2^n}^{2^{n+1}} dx\int_1^{x}\int_1^{x}\frac{dzdy}{\left(x+y\ln x+z\ln^2 x\right)^3}\\
&\leq C\sum_{n=0}^{\infty}\int_1^{2^{n+1}}\int_1^{2^{n+1}}\frac{2^n\, dzdy}{\left(2^{n}+ny+n^2z\right)^3}\\
&\leq C\sum_{n=0}^{\infty}\frac{1}{n^2}<\infty.
\end{align*}
This proves the second part of the theorem. 

Next we  consider Case (b). Since $\alpha_3\leq\alpha_1$, for $|\xi_3|\geq 1$ and $\|\xi\|\geq e$ we have
$$
|\xi_1|^{\alpha_1}+|\xi_2|^{\alpha_1}(\ln\|\xi\|)^{\alpha_1}+|\xi_3|^{\alpha_3}\leq |\xi_1|^{\alpha_1}+|\xi_2|^{\alpha_1}(\ln\|\xi\|)^{\alpha_1}+|\xi_3|^{\alpha_1}(\ln\|\xi\|)^{2\alpha_1}.
$$
Hence, if (\ref{eq:intPsi2}) does not hold in the case (d), then it does not hold in (b) either. Therefore, in what follows we may assume that $\alpha_1\geq3/2$. Since $\alpha_1=\alpha_2=2$, the initial condition $2-1/\alpha_1-1/\alpha_2-1/\alpha_3>0$ translates into $2-2/\alpha_1-1/\alpha_3>0$. By (\ref{eq:Hausdorff2}) and Theorem~\ref{th:Jbeta}, the Hausdorff dimension of $M_2$ is strictly positive (and hence $M_2\neq\emptyset$) if $2-2/\alpha_1-1/\alpha_3>0$, so it is enough to show that (\ref{eq:intPsi2}) does not hold if $2-2/\alpha_1-1/\alpha_3\leq0$. Under this assumption, we need to prove that
$$
\int_1^{\infty}\int_1^{\infty}\int_1^{\infty}\frac{dxdydz}{\left(x^{\alpha_1}+y^{\alpha_1}\ln^{\alpha_1}\sqrt{x^2+y^2+z^2}+z^{\alpha_3}\right)^2}=\infty.
$$
The integral above can be estimated from below by (up to a constant factor),
\begin{align*}
\int_2^{\infty}\int_1^z\int_1^z\frac{ dxdydz}{\left(x^{\alpha_1}+y^{\alpha_1}\ln^{\alpha_1}z+z^{\alpha_3}\right)^2}
&=\sum_{n=1}^{\infty}\int_{2^n}^{2^{n+1}}\int_1^{{z}}\int_1^{z}\frac{dxdydz}{\left(x^{\alpha_1}+y^{\alpha_1}\ln^{\alpha_1}z+z^{\alpha_3}\right)^2}\\
&\geq \sum_{n=1}^{\infty}\int_1^{2^{n}}\int_1^{2^n}\frac{C2^ndxdy}{\left(x+ny+2^{n\alpha_3/\alpha_1}\right)^{2\alpha_1}}.
\end{align*}
Since $\alpha_1\geq3/2$, the last term above is equal to
\begin{equation}\label{eq:SerInt}
C\sum_{n=1}^{\infty}2^n\int_1^{2^{n}}\left[\left(1+ny+2^{n\alpha_3/\alpha_1}\right)^{1-2\alpha_1}
-\left(2^n+ny+2^{n\alpha_3/\alpha_1}\right)^{1-2\alpha_1}\right]dy.
\end{equation}
It is not hard to verify that the last series diverges if $2-2/\alpha_1-1/\alpha_3\leq0$.
This proves the theorem in Case (b).

Finally, we consider Case (c). As in the previous part, we may assume that $\alpha_1\geq3/2$. 
Also, the initial condition $2-1/\alpha_1-1/\alpha_2-1/\alpha_3>0$ becomes 
$2-1/\alpha_1-2/\alpha_2>0$, and by (\ref{eq:Hausdorff2}) and Theorem~\ref{th:Jbeta}, it is 
enough to show that (\ref{eq:intPsi2}) does not hold if $2-1/\alpha_1-2/\alpha_2\leq0$. 
Furthermore, since $\alpha_2\leq\alpha_1$, we may assume that $2-1/\alpha_1-1/\alpha_2>0$. 
Indeed, if $2-1/\alpha_1-1/\alpha_2\leq0$, then the integral in (\ref{eq:intPsi2}) is infinite 
in Case (b), which implies the same for Case (c).  Furthermore it is enough to consider the 
case $\alpha_2<\alpha_1$, since for $\alpha_2=\alpha_1$ the theorem in Case (c) follows 
from that for Case (b), too. We have
\begin{align*}
\int_{\RR^3}\left(\frac{1}{1+\Psi(\xi)}\right)^2d\xi
&\geq \int_2^{\infty}\int_1^{y}\int_1^{y}\frac{Cdxdzdy}{(x^{\alpha_1}+y^{\alpha_2}+z^{\alpha_2}\ln^{\alpha_2}y)^2}\\
&\geq \sum_{n=1}^{\infty}\int_1^{2^{n}}\int_1^{2^n}\frac{C2^ndxdz}{\left(x+2^{n\alpha_2/\alpha_1}+n^{\alpha_2/\alpha_1}z^{\alpha_2/\alpha_1}\right)^{2\alpha_1}}.
\end{align*}
We can verify that the last series is infinite whenever $2-1/\alpha_1-2/\alpha_2\leq0$. This completes the proof of the theorem.
\end{pr}

\bigskip



\bigskip


\begin{thebibliography}{HD}

\bibitem{B96} 
J. Bertoin,  \textit{L\'evy Processes},
    Cambridge Tracts in Mathematics, Cambridge, 1996.
    
\bibitem{DEK1}
A. Dvoretzky, P. Erd\"os, S. Kakutani,
\newblock {\em Double points of paths of Brownian motion in $n$-space},
\newblock  Acta Sci. Math. 12 (1950), 75--81.


\bibitem{DEK2}
A. Dvoretzky, P. Erd\"os, S. Kakutani,
\newblock {\em Multiple points of paths of Brownian motion in the plane},
\newblock Bull. Res. Council Israel, Sect. F 3 (1954), 364--371.


\bibitem{DEKT} A. Dvoretzky, P. Erd\"os, S. Kakutani, S.J. Taylor,
\newblock {\em Triple points of Brownian motion in 3-space},
\newblock Proc. Cambridge Philos. Soc. 53 (1957), 856--862.


\bibitem{E} S.N. Evans,
\newblock {\em Multiple points in the sample paths of a L\'evy process},
\newblock  Probab. Theory Relat. Fields 76 (1987), 359--367.


\bibitem{FS} P.J. Fitzsimmons, T.S. Salisbury,
\newblock {\em Capacity and energy for multiparameter Markov processes},
\newblock  Ann. Inst. H. Poincar\'e Probab. Statist. 25 (1989) 325--350.


\bibitem{H} J. Hawkes,
\newblock {\em Multiple points for symmetric L\'evy processes},
\newblock Math. Proc. Cambridge Philos. Soc. 83 (1978), 83--90.


\bibitem{He} W.J. Hendricks,
\newblock {\em Multiple points for transient symmetric L\'evy processes},
\newblock Z. Wahrsch. verw. Gebiete 49 (1979), 13--21.



\bibitem{K03}
D. Khoshnevisan, {\em Intersections of Brownian motions},
\newblock Expo. Math. 21 (2003),  97--114.

\bibitem{KX4} D. Khoshnevisan, Y. Xiao,
\newblock {\em Level sets of additive L\'evy processes},
\newblock  Ann. Probab. 30 (2002), 62--100.


\bibitem{KX5} D. Khoshnevisan, Y. Xiao,
\newblock {\em Weak unimodality of finite measures, and an application
to potential theory of additive L\'evy processes},
\newblock  Proc. Amer. Math. Soc., 131 (2003), 2611--2616.


\bibitem{KX1} D. Khoshnevisan, Y. Xiao,
\newblock {\em Additive L\'evy processes: capacity and Hausdorff dimension},
\newblock  In: Proc. of Inter. Conf. of Fractal Geometry and Stochastics III.,
Progr. Probab. 57 (2004), 62--100.


\bibitem{KX2} D. Khoshnevisan, Y. Xiao,
\newblock {\em L\'evy processes: capacity and Hausdorff dimension},
\newblock  Ann. Probab. 33 (2005), 841--878.


\bibitem{KX3} D. Khoshnevisan, Y. Xiao,
\newblock {\em Harmonic analysis of additive L\'evy processes},
\newblock  Probab. Theory Relat. Fields 145 (2009), 459--515.

\bibitem{LRS} J.-F. Le Gall, J.S. Rosen, N.-R. Shieh,
\newblock {\em Multiple points of L\'evy processes},
\newblock  Ann. Probab. 17 (1989), 503--515.

\bibitem{MS} M.M. Meerschaert, H.-P. Scheffler,
\newblock {\em Limit Distributions for Sums of Independent Random Vectors},
\newblock  John Wiley, New York, 2001.


\bibitem{MX} M.M. Meerschaert, Y. Xiao,
\newblock {\em Dimension results for sample paths of operator stable L\'evy processes},
\newblock  Stoch. Process. Appl. 115 (2005),  55--75.

\bibitem{S99}
K. Sato,  {\em L\'evy Processes and Infinitely Divisible
Distributions}.  Cambridge Univ. Press, Cambridge, 1999.


\bibitem{S3} N.-R. Shieh,
\newblock {\em Multiple points of dilation-stable L\'evy processes},
\newblock  Ann. Probab. 26 (1998),  1341--1355.

\bibitem{X04}
Y. Xiao, {\em Random fractals and Markov processes}, In: Fractal Geometry and Applications:
A Jubilee of Benoit Mandelbrot, (Michel L. Lapidus and Machiel van Frankenhuijsen, editors),
pp. 261--338, American Mathematical Society, 2004.


\end{thebibliography}
\end{document}